\documentclass[ijoc,nonblindrev]{informs3_2} 

\usepackage{multirow, makecell}
\usepackage{booktabs} 
\usepackage[linesnumbered,ruled,vlined]{algorithm2e}

\usepackage{natbib}
\bibpunct[, ]{(}{)}{,}{a}{}{,}%
\def\bibfont{\small}%
\TheoremsNumberedThrough     

\EquationsNumberedThrough    

\usepackage{subfig, float}
\usepackage{multirow, makecell}
\usepackage[dvipsnames]{xcolor}
\usepackage[prependcaption,textsize=small]{todonotes}

\newcommand{\A}{\mathcal A}
\newcommand{\X}{\mathcal X}
\newcommand{\I}{\mathcal I}
\newcommand{\J}{\mathcal J}
\newcommand{\B}{\mathcal B}
\newcommand{\Obj}{\mathfrak O}
\newcommand{\U}{\mathfrak U}
\newcommand{\R}{\mathcal R}
\newcommand{\C}{\mathcal C}
\newcommand{\K}{\mathcal K}
\newcommand{\Q}{\mathcal Q}

\newcommand{\V}{\mathcal V}
\newcommand{\E}{\mathcal E}
\newcommand{\N}{\mathcal N}

\newcommand{\Uc}{\mathcal U}
\newcommand{\G}{\mathcal G}
\newcommand{\T}{\mathcal T}
\begin{document}
	
	
	\RUNAUTHOR{Byeon, Van Hentenryck}
	
	\RUNTITLE{Benders Subproblem Decomposition for Bilevel Problems with Convex Follower}
	
	\TITLE{Benders Subproblem Decomposition \\ for Bilevel Problems with Convex Follower}
	
	\ARTICLEAUTHORS{%
		\AUTHOR{Geunyeong Byeon}
		\AFF{School of Computing, Informatics, and Decision Systems Engineering, Arizona State University \EMAIL{} \URL{}}
		\AUTHOR{Pascal Van Hentenryck}
				\AFF{H. Milton Stewart School of Industrial and Systems  Engineering,  Georgia  Institute  of  Technology \EMAIL{} \URL{}}
	} 

	\ABSTRACT{%
	Bilevel optimization formulates hierarchical decision-making processes that arise in many real-world applications such as in pricing, network design, and infrastructure defense planning. In this paper, we consider a class of bilevel optimization problems where the upper level problem features some integer variables while the lower level problem enjoys strong duality. We propose a dedicated Benders decomposition method for solving this class of bilevel problems, which decomposes the Benders subproblem into two more tractable, sequentially solvable problems that can be interpreted as the upper and the lower level problems. We show that the Benders subproblem decomposition carries over to an interesting extension of bilevel problems, which connects the upper level solution with the lower level dual solution, and discuss some special cases of bilevel problems that allow sequence-independent subproblem decomposition. Several novel schemes for generating numerically stable cuts, finding a good incumbent solution, and accelerating the search tree are discussed. A computational study demonstrates the computational benefits of the proposed method over a state-of-the-art bilevel-tailored branch-and-cut method, a commercial solver, and the standard Benders method on standard test cases and the motivating applications in sequential energy markets. 
	}%
	
	
	
	\maketitle

\section{Introduction}

 A variety of real-world applications involves multiple decision makers. These decision makers (agents) may have an implicit \emph{hierarchy} in the sense that the decision made by an agent at a certain level of the hierarchy precedes and affects the decisions of agents at lower levels that, in turn, affect the outcomes of the decisions at the higher levels. \emph{Hierarchical optimization} models optimization problems that involve the hierarchical decision-making process of multiple agents.


Bilevel optimization is a subclass of hierarchical optimization with a two-level decision hierarchy, the upper- and lower-level of which is often referred to as a leader and a follower, respectively. In these problems, it is assumed that the leader can anticipate how the follower would respond to her decision. The objective of these problems is, thus, to find an optimal decision for the leader by solving an optimization problem that embeds the follower problem as a constraint (see, e.g., \cite{dempe2002foundations} for more details). In this paper, we consider a family of bilevel optimization problems in which the leader problem is modeled as a Mixed-Integer Second-Order Cone Programming (MISOCP) problem and the follower problem is modeled as a Second-Order Cone Programming (SOCP) problem. We name this problem class Bilevel Second-Order Cone Programming (BSOCP)\footnote{For clarification, we denote a bilevel optimization problem as `mixed-integer' only when both the leader and follower are allowed to have integer variables. The class of bilevel problems of interest assumes continuous follower variables, so we denote the class as BSOCP, even though it is allowed to have integer variables in the leader problem.}.

BSOCPs arise in many applications across various sectors including economics, energy infrastructure, and defense. For instance, a special class of BSOCP with only linear constraints, referred to as Bilevel Linear Programming (BLP), models various network planning/design problems with autonomous agents, e.g., the optimal zonal configuration problem in zonal-pricing electricity markets \citep{grimm2019optimal, ambrosius2018endogenous, kleinert2019global}, an urban traffic network design problem \citep{fontaine2014benders}, and facility location problems for logistics distribution center \citep{sun2008bi}. BLP can also be used to model the evasive flow capturing problem discussed by \cite{arslan2018exact} which has applications in transportation, revenue management, and security management. In addition, unit scheduling problems under sequentially cleared markets can be modeled with BSOCP in its extended form where an additional constraint stating the impact of the \emph{follower's dual solution} on the leader problem is added to the leader problem. Such constraints may be desirable in some sequential market environments where the follower's dual solution settles the prices of commodities that are used by the leader; see, e.g., the unit-commitment problem for interdependent natural gas and electricity markets studied by \cite{byeon2019unit}.

Despite the wide applicability of bilevel optimization to problems with multiple decision makers, the nonconvexity brought by the optimality requirement of the follower makes even the simplest subclass of bilevel problems, i.e., linear leader and follower problems, NP-hard \citep{jeroslow1985polynomial}. This inherent complexity of bilevel optimization explains why the design of tailored algorithms for bilevel optimization problems with integer variables has only a relatively short history \citep{denegre2009branch,xu2014exact,fischetti2016intersection,fischetti2017new,lozano2017value,kleinert2020outer}. 
These dedicated methods adapt branch-and-bound and/or cutting-plane approaches to the bilevel context. However, most of them focus on the case of \emph{linear} or \emph{convex quadratic} constraints. Moreover, only the work of \cite{lozano2017value} considers continuous non-linear constraints, but it requires all leader variables to be integer-valued for finite termination.

Aside from these bilevel-tailored branch-and-bound and/or cutting-plane methods, an alternative approach for solving BSOCP is to reformulate the bilevel optimization problem into a single-level optimization problem. \cite{cao2006capacitated, garces2009bilevel,fontaine2014benders}, and \cite{zare2019note} reformulated a BLP instance as a single-level Mixed-Integer Linear Programming (MILP) problem by replacing the lower level with its optimality conditions. The same technique can be applied to BSOCP, leading to a single-level MISOCP reformulation. The resulting MISOCP problem, however, is complex since it intertwines the leader problem and the follower's primal and dual problems. For large instances, the complexity of the MISOCP formulation often raises significant challenges for commercial solvers.

Benders decomposition is a solution technique that has been widely used for solving large-scale MILP and MISOCP problems. Instead of dealing with all the variables and constraints of a complex MISOCP problem simultaneously, Benders decomposition relaxes the inner-continuous problem and iteratively discovers the shape and the domain of the relaxed problem with a set of linear inequalities called optimality and feasibility cuts. At each iteration, the cuts are generated by solving a Benders subproblem. However, the complexity of the MISOCP formulation of BSOCP also complicates the Benders subproblem, which often exhibits numerical difficulties and requires significant computational resources.

To address these challenges, we propose a dedicated Benders decomposition for BSOCP 
where the complex Benders subproblem is itself decomposed into two more tractable, sequentially solvable problems that are closely related to the leader and the follower problems. Moreover,
to address applications where the dual variables of the follower problem have no natural bounds, we propose a new family of cuts that merges no-good and Benders cuts, removing the need for these bounds and reducing both the computational burden and the numerical issues. Since this novel decomposition is embedded into a branch and cut algorithm, we also propose two novel techniques to speed up the solution space. First, we propose a new branching scheme that targets the optimality gap between the follower objective and its guess in the leader subproblem. Second, we propose a new
method for finding a good incumbent solution in a preprocessing step, which combines the bilevel-tailored branching scheme and a heuristic local cut. We also show that the proposed decomposition applies to applications where the leader problem features constraints on the dual variables of the follower problem, which is the case in our motivating case study in sequential market clearing for electricity and gas networks. We also identify cases where the two subproblems can be solved independently (instead of sequentially) and propose acceleration schemes that improve the performance of the overall algorithm in this setting. 

To the best of our knowledge, tailored algorithms for BSOCP, especially for the case where the follower problem is a general SOCP, have not been discussed or their computational studies have focused only on linear cases. Due to the lack of dedicated algorithms and available code packages for solving the BSOCP instances, we benchmarked the proposed method against a single-level MISOCP reformulation which will be described in Section \ref{sec:MIP}. In addition, to demonstrate the potential benefits of our proposed approach, we also conducted extensive experiments on BLP instances, a special case of BSOCP, which were constructed by relaxing the follower integrality condition from publicly available MIBLP test sets. The performance of our proposed approach is then compared to a state-of-the-art bilevel-tailored branch-and-cut algorithm that can be used to solve MIBLP \citep{fischetti2017new}.

The main contributions of this paper can be summarized as follows.
\begin{itemize}
    \item The paper proposes a new decomposition technique for BSOCP, which allows for easy implementation and an intuitive interpretation of Benders cuts.
    \item The paper proposes a new family of hybrid cuts that combine no-good and Benders cuts to eliminate the need for bounds on the dual variables of the follower problem.  
    \item The paper proposes a new branching scheme for BSOCP that targets the optimality gap between the follower objective and its guess in the leader problem, as well as a new method for finding a high-quality solution before the branch and cut exploration. 
    \item The paper introduces an interesting extension of BSOCP that captures important real-world problems where the leader is affected by the follower's dual solution.  It is shown that the decomposition technique carries over to this extension. In addition, the paper identifies special cases of BSOCP that allow for a sequence-independent decomposition. The paper also presents some accelerating schemes to further reduce the computational burden in this case.
    \item The paper reports a computational study that demonstrates significant performance improvement of the proposed method and the accelerating schemes over a commercial solver and the standard Benders method. Extensive experiments on BLP are also given: they highlight the benefits of the proposed method and its complementarity with a state-of-the-art bilevel-tailored branch-and-cut algorithm \citep{fischetti2017new}.
\end{itemize}

The rest of the paper is organized as follows: Section \ref{sec:assum} formally defines BSOCP, as well as the assumptions of the paper and their justifications. Section \ref{sec:review} discusses previous work and Section \ref{sec:MIP} presents the MISOCP formulation of BSOCP. Sections \ref{sec:decom} and \ref{sec:stable} propose the dedicated Benders method for BSOCP, and its numerically stable variant. Section \ref{sec:heuristic} proposes a heuristic method for finding an incumbent solution in a preprocessing step. Section \ref{sec:ext} discusses an interesting extension of BSOCP that incorporates additional constraints on the follower's dual variables in the leader problem. It also identifies some special cases of BSOCP that allow stronger algorithmic results, as well as some accelerating schemes for the dedicated Benders method. The computational performance of the proposed method is demonstrated in Section \ref{sec:comp}. Section \ref{sec:conclusion} concludes the paper.

\subsection{Bilevel Secone-Order Cone Programming (BSOCP) and Assumptions}\label{sec:assum}
A BSOCP problem is formally defined as follows:
\begin{subequations}
\begin{alignat}{4}
 \quad &\min_{x,y} \qquad && c_x^T x + c_y^T y\label{prob:bl:obj}\\
& \mbox{ s.t.}    && G_{x} x + G_y y \ge h, \label{prob:bl:upper:y}\\
& &&  x \in \X:=\{x \in \K_x: x_i \in [\underline{x}_i, \overline{x}_i]_{
\mathbb Z}, \ \forall i \in \I\},\label{prob:bl:upper:x}\\
& && y \in \arg \min_{y \in \K_y} \{d^T y: Ax + By \ge b \}, \label{prob:bl:lower}
\end{alignat}\label{prob:bl}
\end{subequations}
\noindent
where $x$ and $y$ respectively represent the $n_x$-dimensional leader and $n_y$-dimensional follower variables. In Equation \eqref{prob:bl:upper:x}, $\underline{x}_i$ and $\overline{x}_i$ respectively denote lower and upper bounds on variable $x_i$, some of which are allowed to be $-\infty$ and $\infty$. $[\underline{x}_i,\overline{x}_i]_{\mathbb Z}$ denotes a set of integer points in the interval, and $\I \subseteq \{1, \cdots, n_x\}$ represents a set of indices of the leader's variables for which the corresponding variable is integer. Each of $\K_x \subseteq \mathbb R^{n_x}$ and $\K_y \subseteq \mathbb R^{n_y}$ is the Cartesian product of a collection of second-order cones and nonnegative orthants, i.e., 
\[\K_x \times \K_y = \K_{n_1} \times \cdots \times \K_{n_l}\]
where each $\K_{n_i} \subseteq \mathbb R^{n_i}$ is either a $n_i$-dimensional second-order cone $\{(u, v) \in \mathbb R^{n_i}: \|u\|_2 \le v\}$ or a $n_i$-dimensional nonnegative orthant $\mathbb R_+^{n_i}$.
$G_x \in \mathbb R^{m_x \times n_x}, G_y \in \mathbb R^{m_x \times n_y}, c_x \in \mathbb R^{n_x}, c_y \in \mathbb R^{n_y}, h \in \mathbb R^{m_x}, d \in \mathbb R^{n_y}, A \in \mathbb R^{m_y \times n_x}, B \in \mathbb R^{m_y \times n_y}$, and $b\in \mathbb R^{m_y}$ are given rational matrices or vectors. 

\begin{remark}
	Note that the leader problem minimizes over $x$ and $y$, which implies a cooperative behavior of the leader and the follower, i.e., when there are multiple lower level optimal solutions for a given upper-level decision $\hat x$, it chooses $\hat y$ that benefits the upper-level the most, among the follower optimal solutions. Bilevel problems with this property are said to be \emph{optimistic}; For more details on this topic, we refer the reader to \citet{colson2005bilevel}.
\end{remark}

Throughout this paper, we assume the following:
 \begin{assumption}
 The dual of the follower problem is feasible when variables $x$ are assigned to a leader decision:
 \label{assum:sd}
 \end{assumption}
The dual of the follower problem for given $\hat{x}$ is
\begin{equation}\max_{\psi \in \mathbb R^{m_y}_{+}} \ (b-A\hat{x})^T\psi: B^T \psi \preceq_{\K_y} d.
\label{prob:follower:dual}
\end{equation}
Note that the dual feasible region is not affected by $\hat{x}$, and thus this assumption implies that the follower problem is bounded from below for any $\hat{x}$. Therefore, for any given $\hat{x}$, strong duality holds between the primal and dual problems \eqref{prob:bl:lower} and \eqref{prob:follower:dual}. This is a reasonable assumption since, otherwise, the follower problem is either unbounded or infeasible for any leader decision.

Let $\J$ denote the set of indices of the leader variables that appear in the follower problem, i.e., $i \in \J$ if and only if the $i$-th column of $A$ is nontrivial (a nonzero vector).
\begin{assumption}
(a) $\J \subseteq \I$, and (b) for each $i \in \J$, $\underline x_i$ and $\overline x_i$ are finite real numbers.
\label{assum:binary_upper_level}
\end{assumption}
This assumption is required by many other state-of-the-art algorithms for mixed-integer bilevel problems, e.g, \citep{xu2014exact,fischetti2016intersection,fischetti2017new,lozano2017value,kleinert2020outer}.

Consider a single-level optimization problem that gives a lower bound to the bilevel program, the so-called \emph{high point problem} (HPP), that is obtained by relaxing the optimality requirement of the follower: \begin{subequations}
\begin{alignat}{4}
	\quad &\min_{x \in \X,  y \in \mathbb \K_y} \qquad && c_x^T x + c_y^T y\\
	& \mbox{ s.t.}    && G_{xy} x + G_y y \ge h_y, \\
	& && Ax + By \ge b.
\end{alignat}\label{prob:hpp}
\end{subequations}
\noindent In order to ensure that Problem \eqref{prob:bl} is neither infeasible nor unbounded, we make the following additional assumptions:
\begin{assumption}
Problem \eqref{prob:hpp} has a bounded feasible region.
\label{assum:lb}
\end{assumption} 
Assumption \ref{assum:lb} guarantees a finite lower bound of Problem \eqref{prob:bl}; this assumption is not too restrictive, because we can add auxiliary variables with penalties to guarantee feasibility and many real-world applications have natural bounds on variables. This assumption holds in energy systems where it is always possible to shed the load, albeit with a high penalty. With these assumptions, we use the following definition throughout this paper:
\begin{definition} 
	A leader decision $\hat x$ is called \emph{bilevel-feasible}, if it satisfies all of the following: 
	\begin{enumerate}
		\item $\hat x$ is feasible to the follower (i.e., Problem \eqref{prob:bl:lower} with $x$ fixed as $\hat x$ is feasible);
		\item there exists an optimal response $\hat y$ of the follower to $\hat x$ that is feasible to the leader, i.e., $G_{xy} \hat x + G_y \hat y \ge h$ and $\hat y \in \mathcal F(\hat x)$, where $\mathcal F(\hat x)$ denotes the set of optimal solutions of the follower problem for given $\hat x$. 
	\end{enumerate}
	When at least one of the above is not met, we call the pair \emph{bilevel-infeasible}.
	\label{def:bilevel_feasible}
\end{definition}

Additionally, based on Assumption \ref{assum:binary_upper_level} (b), we assume, w.l.o.g, that $\underline x_i$ and $\overline x_i$ are integers for $i \in \J$. Note that an integer variable $x_i$ with finite integral upper and lower bounds can be replaced by a set of auxiliary binary variables $\{z_{ij}\}_{j = 1, \cdots, k_i}$, where $k_i := \lfloor\log(\overline {x}_i - \underline {x}_i)\rfloor+1$, as follows: $x_i = \sum_{j = 1}^{k_i} 2^{j-1}z_{ij} + \underline x_i$. Therefore, w.l.o.g, we assume $x_i \in \{0,1\}, \forall i \in \J$.

	\section{Literature Review}\label{sec:review}



A widely-studied special class of BSOCP is BLP where both $\K_x$ and $\K_y$ are $n_x$- and $n_y$-dimensional nonnegative orthants, respectively. Taking advantage of the strong duality in the lower level problem, the common solution approach for BLP is to reformulate the bilevel problem into a single-level MILP problem and to solve the MILP problem
via off-the-shelf solvers. There are two widely-used reformulation
schemes: (1) a Karush-Kuhn-Tucker (KKT) condition approach, and (2) a
strong duality approach. The former replaces the lower level problem
by the KKT conditions and linearizes the nonlinear complementary
slackness condition by introducing additional binary variables and
logic-based constraints (see, e.g., \cite{labbe1998bilevel}). However,
due to the large number of binary variables and constraints that
should be introduced for the linearization, this approach does not
scale well and is not adequate for solving large instances. The
later method, on the other hand, replaces the complementary slackness
condition with the reversed weak duality inequality to ensure that the
primal and dual objective values of the lower level are the
same. Then, the bilinear terms in the reversed weak duality are
linearized using the McCormick relaxation
\citep{cao2006capacitated,garces2009bilevel, fontaine2014benders} or
some problem-specific properties \citep{arslan2018exact}. Recently,
\cite{zare2019note} have compared these two schemes and
have shown that the latter approach significantly outperforms
 the former approach for many classes of instances.

For large-scale problems, however, solving the resultant MILP is still
challenging since it entangles the leader problem and the follower
primal and dual problems. Accordingly, some problem-specific and
generic decomposition/separation techniques for solving the associated
MILP have been proposed. \cite{grimm2019optimal} proposed a generalized
Benders algorithm that uses a special structure of the given tri-level
problem (which has an equivalent BLP counterpart) and
\cite{arslan2018exact} developed a branch-and-cut approach for a certain class
of BLP, named the Evasive Flow Capturing Problem. For general
approaches, \cite{saharidis2009resolution} proposed a decomposition
algorithm which, at every iteration, fixes the integer variables at
some values, reformulates the resultant bilevel linear subproblem into
a MILP problem using the KKT scheme, solves the MILP problem to construct the associated LP problem with its active constraint set, solves the LP problem to obtain the dual information, and adds a
cut. Since this approach reformulates the bilevel linear subproblem as
a MILP problem using the KKT scheme at every iteration, its application
to large-scale problems would be computationally expensive. The most
relevant work is by \cite{fontaine2014benders} who applied the Benders
decomposition to the MILP formulation obtained by the strong duality
scheme. It proposed an acceleration scheme for obtaining an optimality Benders cut which sequentially solves three smaller problems: (a) the follower problem, (b) the leader
problem, and (c) the follower dual-related problem to obtain optimality
cut. 

Another line of research has developed bilevel-tailored branch-and-bound and/or cutting plane methods for Mixed-Integer Bilevel Linear Programming (BMILP), which subsumes BLP, where some of the leader and follower
variables are allowed to be integer-valued. \cite{xu2014exact} proposed a branch-and-bound approach which features a bilevel-tailored design of the relaxation problem and the branch-and-bound rules. 
Separately, in the spirit of cutting-plane approach, \cite{denegre2009branch} and \cite{fischetti2016intersection} proposed valid cut generation schemes for Integer Bilevel Programming (IBLP) and MIBLP respectively, which were further improved by \cite{fischetti2017new}. The Branch-and-Cut (B\&C) approach proposed by \cite{fischetti2017new} features intersection cuts along with two additional acceleration schemes: (i) locally valid cuts and (ii) a preprocessing rule that allows for the pre-determination of some follower solutions, which enhances the algorithmic performance significantly. Another
B\&C algorithm for MIBLP was proposed by \cite{caramia2015enhanced}, which solves a BLP
problem for generating cuts.
Although, these generic methods solved quite large instances of MIBLP, they focused on the case of \emph{linear} follower constraints. Only the work of \cite{kleinert2020outer} and \cite{lozano2017value} considered a non-linear follower problem: \cite{kleinert2020outer} developed an outer-approximation-based cutting-plane method for solving a special class of BSOCP in which some of its constraints are allowed to be convex quadratic and \cite{lozano2017value} proposed a sampling-based cutting plane method for solving a class of mixed-integer nonlinear bilevel programming, where all the leader variables are assumed to be integers. 

From a computational standpoint, this paper features three main differences from the existing literature: (i) it proposes a tailored solution approach for solving BSOCP problems; (ii) it presents a new family of cuts that combines no-good and Benders cuts, removes the need for bounds of the follower dual variables, and leads to better performance, numerical stability, and ease of implementation; (iii) it develops a heuristic method for finding a good incumbent solution of MIBSOCP. From a modeling standpoint, the paper shows that the decomposition carries over to the interesting case where the leader problem feature constraints on the follower dual variables. Additional modeling and computational results are also presented in this more general setting.
	\section{The MISOCP Reformulation}\label{sec:MIP}
In this section, we reformulate Problem \eqref{prob:bl} as a single-level MISOCP problem using the strong duality approach. Note that, using Assumptions \ref{assum:sd}, Problem \eqref{prob:bl} can be expressed as follows:\\
\noindent\begin{minipage}{0.4\textwidth}
    \centering
    \begin{subequations}
        \begin{alignat}{4}
         \min_{x \in \mathcal X, t \in \mathbb R} \ & c_x^T x  +  t\\
         \mbox{s.t.} \ &t \ge f(x), \label{prob:bl:1st:constr} 
        \end{alignat}\label{prob:bl:1st}
        \end{subequations}  
  \vspace{2.5cm}      
\end{minipage}
\begin{minipage}{0.5\textwidth}
    \centering
    \begin{subequations}
\begin{alignat}{3}
 f(x) := \ &\min_{y \in \K_y,\psi \ge 0} \quad && c_y^T y & \label{eq:bl:inner:obj}\\
& \qquad \mbox{ s.t.}    && G_{xy} x + G_y y \ge h_y, \label{eq:bl:inner:1st:y}\\
& && Ax + By \ge b,\label{eq:bl:inner:2nd:primal}\\
& && B^T \psi \preceq_{\K_y} d,\label{eq:bl:inner:2nd:dual}\\
& && d^T y \le \psi^T  (b - A x).\label{eq:bl:inner:sd}
\end{alignat}\label{prob:bl:inner:bilinear}
\end{subequations}
\end{minipage}\\
\noindent Constraints \eqref{eq:bl:inner:2nd:primal} and
\eqref{eq:bl:inner:2nd:dual} respectively ensure primal and dual
feasibility of the lower level problem, Constraint
\eqref{eq:bl:inner:sd} ensures strong duality in the lower
level. Thus, for any $x \in \mathbb R^{n_1}$, a feasible $y$ to Constraints \eqref{eq:bl:inner:2nd:primal}-\eqref{eq:bl:inner:sd} is an optimal solution of the lower level
problem for the given $x$. Accordingly, Constraint
\eqref{eq:bl:inner:1st:y} models how the lower level reaction affects
the upper level feasible region. 

Problem \eqref{prob:bl:inner:bilinear} contains a bilinear term, $\psi^T
A x = \sum_{i=1}^{m_y} \sum_{j=1}^{n_x} A_{ij}\psi_i x_j$, in Constraint \eqref{eq:bl:inner:sd}. Note that, due to Assumption \ref{assum:binary_upper_level}, each non-trivial bilinear term $A_{ij}\psi_i x_j$ is a multiplication of some nonnegative continuous variable $\psi_i$ and a binary variable $x_j$.
Assuming that $\psi$ has an
upper bound of $\overline \psi$\footnote{The case where a reasonable $\overline \psi$ is not available is discussed in Section \ref{sec:stable}}, each of the nonlinear terms can be linearized. First,
introduce an additional vector of nonnegative variables $\mu \in \mathbb {R}_+^{m_y|\J|}$ and constraints $
    \mu_{(i-1)|\J| + j} = \psi_i x_j, \forall i=1, \cdots, m_y, \ j \in \J$ 
to represent $\psi^T A x$ as $\mu^T a$, where $a$ is a vector obtained by concatenating each rows of $A$. Then, for each $i=1,\cdots,m_y$ and $j \in \J$, use a McCormick
transformation to replace the additional constraint by a
set of linear constraints of the form:
$-\psi_{i} + \mu_{(i-1)|\J| + j}  \ge \overline{\psi}_i x_{j}-\overline{\psi}_i, \
\mu_{(i-1)|\J| + j} \le \overline{\psi}_i x_{j}, \
-\psi_{i} +\mu_{(i-1)|\J| + j} \le 0.
\label{eq:Mc0}
$
We represent this set of equations for all $i=1,\cdots,m_y$ and $j \in \J$ as
\begin{equation}
K_\psi \psi + K_\mu \mu \ge k + K_x x,
\label{eq:Mc}
\end{equation}
for some matrices $K_\psi, K_\mu, K_x,$ and some vector $k$ of appropriate dimensions. Then, $f(x)$ can be
obtained by solving the following problem:
\begin{subequations}
\begin{alignat}{4}
 &\min_{(y, \psi, \mu)^T \in \K_y \times \mathbb R_+^{m_y} \times  \mathbb {R}_+^{m_y|\J|}} \qquad && c_y^T y   \\
& \mbox{ s.t.}  &&  G_{y} y \ge h_y - G_{xy} x, \label{eq:inner:y}\\
& && By \ge b - A x, \label{eq:inner:2nd:primal}\\
& && -\psi^T B \succeq_{\K_y} -d^T,\label{eq:inner:2nd:dual}\\
& && -d^T y + \psi^T b - \mu^T a \ge 0,\label{eq:inner:2nd:sd}\\
& && K_\psi \psi +K_\mu \mu \ge k + K_x  x,\label{eq:inner:2nd:Mc}.
\end{alignat}\label{prob:inner}    
\end{subequations}
\noindent In the following, Problem (MISOCP) denotes the resulting MISOCP problem, i.e., Problem \eqref{prob:bl:1st} where $f(x)$ is defined by Problem \eqref{prob:inner}.

	\section{A Dedicated Benders Decomposition Method for BSOCP}
\label{sec:decom}
This section discusses a dedicated solution method for Problem (MISOCP), which builds upon the Benders decomposition method---a solution technique that has been widely used for solving large-scale MILP and MISOCP problems. Benders Decomposition (BD) is defined by a Relaxed Master Problem (RMP) and a Benders SubProblem (BSP).
Initially, the RMP corresponds to Problem \eqref{prob:bl:1st} with Constraint \eqref{prob:bl:1st:constr} relaxed:
\begin{equation}
	\begin{alignedat}{2}
	   &\min_{x \in \X} \quad && c_x^T  x + t\\
	   &\mbox{ s.t.}  && t \in \mathbb R.
	\end{alignedat}\label{prob:M0}
\end{equation}

\begin{figure}[!tb]
	\centering
	  \subfloat[A feasibility cut]{\includegraphics[width=0.3\textwidth]{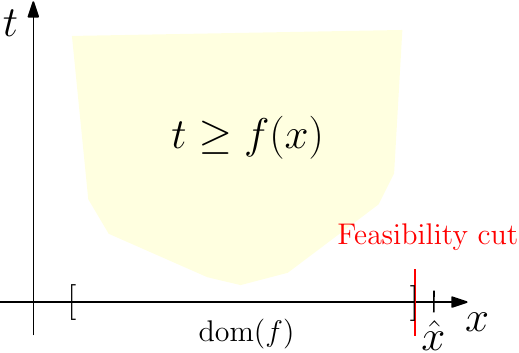}}
	  \subfloat[An optimality cut]{\includegraphics[width=0.35\textwidth]{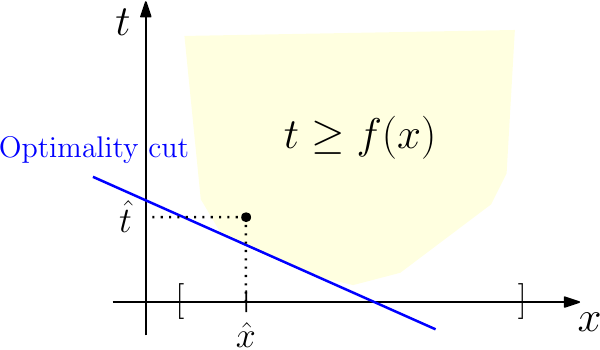}}
	\caption{Illustration of Benders Cuts\label{fig:benders}}
	\end{figure}

At each iteration, BD generates a guess $(\hat x, \hat t)$ by solving the RMP and then checks whether $(\hat x, \hat{t})$ violates the relaxed constraint or not by solving the BSP, which is defined by the dual of $f(\hat x)$. There are three possible cases: (i) $f(\hat x) = \infty$ (i.e., $\hat x$ is out of the domain of $f$); (ii) $f(\hat x)< \infty$ but $(\hat x, \hat t)$ violates Constraint \eqref{prob:bl:1st:constr}; (iii) $f(\hat x)< \infty$ and $(\hat x, \hat t)$ satisfies Constraint \eqref{prob:bl:1st:constr}. For infeasible cases (i) and (ii), BD respectively generates feasibility and optimality cuts using an unbounded ray and an optimal solution of the BSP to cut off the current guess, which is illustrated in Figure \ref{fig:benders}. BD repeats this procedure until it reaches a predetermined gap or encounter case (iii). The detailed idea behind the Benders cuts can be found in many optimization textbooks (e.g., \cite{wolsey1998integer,conforti2014integer,junger200950}).

Specifically, for a guess $\hat x$, the BSP is defined by the dual of Problem
\eqref{prob:inner}. We slightly abuse notation and let $u_y$, $\psi$,
$y$, $w$, and $v$ represent the dual variables associated with
Constraints \eqref{eq:inner:y}, \eqref{eq:inner:2nd:primal},
\eqref{eq:inner:2nd:dual}, \eqref{eq:inner:2nd:sd}, and
\eqref{eq:inner:2nd:Mc} respectively. Then the dual of Problem
\eqref{prob:inner} can be expressed as follows:
\begin{subequations}
\begin{alignat}{4}
 &\max \qquad && \psi^T (b - A \hat x) + u_y^T (h_y - G_{xy} \hat x) - \left[d^T y - v^T (k + K_x \hat x)\right]\\
& \mbox{ s.t.}  &&  By  - K_\psi^T v \ge bw, \label{eq:inner:dual:primal1}\\
& && B^T \psi + G_y^T u_y  \preceq_{\K_y} dw + c_y,\label{eq:inner:dual:dual}\\
& && K_\mu^T v \le a w,\label{eq:inner:dual:primal2}\\
& && \psi \ge 0,u_y \ge 0,w \ge 0,y \in \K_y,v \ge 0.
\end{alignat}\label{prob:inner:dual}
\end{subequations}
\noindent
 Note that, without loss of generality, we can assume that Problem \eqref{prob:inner:dual} is feasible, since otherwise, Problem \eqref{prob:bl:inner:bilinear} is infeasible for any $\hat x \in \X$ due to Assumption \ref{assum:lb}, and thus we can conclude that Problem \eqref{prob:bl} is infeasible. 

Unfortunately, for large-scale bilevel problems, Problem
\eqref{prob:inner:dual} is highly complex since it has
primal-related (e.g., \eqref{eq:inner:dual:primal1} and
\eqref{eq:inner:dual:primal2}) and dual-related (e.g.,
\eqref{eq:inner:dual:dual}) constraints for Problem \eqref{prob:bl:lower} which are linked by variable $w$. In this section, we show that Problem
  \eqref{prob:inner:dual} does not need to be solved as a whole. Rather, the
Benders cuts of Problem (MISOCP) can be obtained by solving two
more tractable problems, i.e.,
a problem associated with the lower-level problem (to be defined as Problem \eqref{prob:s1}) and a problem related to the upper level problem (to be defined as Problem \eqref{prob:s2}).

\begin{theorem}
  Problem \eqref{prob:inner:dual} can be solved by solving two more tractable problems sequentially, i.e., solve the following problems \\
  \begin{minipage}{0.45\textwidth}
	\begin{subequations}
	\begin{alignat}{2}
	\min \ &  d^T y  - v^T (k+K_x\hat x)  \\
    \mbox{ s.t.}  \ & By  - K_\psi^T v \ge b, \label{eq:s1:1}\\
     & K_\mu^T v \le a,\label{eq:s1:2}\\
     & y \in \K_y,v \ge 0,
	\end{alignat}\label{prob:s1}
	\end{subequations}
\end{minipage}
\begin{minipage}{0.55\textwidth}
	\begin{subequations}
	\begin{alignat}{2}
	\max \ &  \psi^T (b - A \hat x) + u_y^T (h_y - G_{xy} \hat x) - \Obj w \\
    \mbox{ s.t.}  \ &B^T \psi + G_y^T u_y  \preceq_{\K_y} dw + c_y, \label{eq:s2:1}\\
     & \psi \ge 0,u_y \ge 0,w \ge 0,
	\end{alignat}\label{prob:s2}
	\end{subequations}
	\end{minipage}
where $\Obj$ denotes the optimal objective value of Problem \eqref{prob:s1}.
\label{theo:bender}
\end{theorem}

\begin{remark} Note that Problem \eqref{prob:s1} has a finite optimum $\Obj$ for any $\hat x$. Consider the dual of Problem \eqref{prob:s1}: 
$\max_{\psi\ge 0,\mu \ge 0} \left\{  b^T \psi - \mu^T a: B^T \psi \preceq_{\K_y} d, K_\psi \psi + K_\mu \mu \ge k + K_x \hat x\right\}$, where $\psi$ and $\mu$ are dual variables associated with Constraints \eqref{eq:s1:1} and \eqref{eq:s1:2} respectively. Note that, due to the McCormick relaxation, it holds that 
	\begin{equation}
		\mbox{Dual of Problem \eqref{prob:s1}} \ge \max_{0 \le \psi \le \bar \psi} \{\psi^T (b - A\hat x) : B^T \psi \preceq_{\K_y} d\},
	\label{prob:follower:dual:bounded}
	\end{equation}
	where the inequality holds in equality for any $\hat x \in \X$, (i.e., when the integrality condition is met). Note that Problem \eqref{prob:follower:dual:bounded} has a nonempty bounded feasible region, the feasibility of which is guaranteed by Assumption \ref{assum:sd} and the boundedness follows from $0 \le \psi \le \bar \psi$. Therefore, Problem \eqref{prob:s1} is always bounded below. Note also that Problem \eqref{prob:s1} is feasible; otherwise, Problem \eqref{prob:s1} is infeasible for any $\hat x \in \mathbb R^{n_x}$, and thus Problem \eqref{prob:follower:dual:bounded} is infeasible or unbounded for any $\hat x \in \X$, which contradicts Assumption \ref{assum:sd}.
	
	Likewise, Assumption \ref{assum:lb} guarantees Problem \eqref{prob:s2} to be feasible. Consider the dual of Problem \eqref{prob:s2}:
\begin{equation}
    \min_{y \in \K_y} \{c_y^T y:By \ge b - A \hat x, \
G_y y \ge h_y - G_{xy} \hat x,\ d^T y \le \Obj\}.
\label{prob:upper:fixed:primal}
\end{equation}
	Note that, if Problem \eqref{prob:s2} is infeasible, Problem \eqref{prob:upper:fixed:primal} is infeasible or unbounded for any $\hat x \in \X$ and $\Obj \in \mathbb R \cup \{\infty\}$, which contradicts Assumption \ref{assum:lb}.
\label{rem:theo1:feas}
\end{remark}

Theorem \ref{theo:bender} implies that Benders cuts can be generated by solving Problem \eqref{prob:s1} (i.e., a lower level-related problem) and Problem \eqref{prob:s2} (i.e., an upper level-related problem) sequentially, and leads to the following corollary.

\begin{corollary}\label{coro:1}
Problem (MISOCP) is equivalent
to the following problem:
{\fontsize{11}{17}\selectfont
	\begin{subequations}
		\begin{alignat}{3}
		\min_{x \in \X} \ & c_x^T x + t && \nonumber\\
		  \mbox{ s.t.} \ & t \ge \hat\psi^T (b- Ax) + \hat u_y^T (h_y -G_{xy} x) - \hat w(d^T \hat y  -\hat v^T(k+K_x x)), &&\forall (\hat \psi, \hat u_y, \hat w,\hat y, \hat v) \in \J_2 \times \J_1, \label{eq:opt_cut}\\
		 & 0 \ge \tilde\psi^T (b- Ax) + \tilde u_y^T (h_y -G_{xy} x) - \tilde w(d^T \hat y  -\hat v^T(k+K_x x)), \ && \forall (\tilde \psi, \tilde u_y, \tilde w, \hat y,\hat v) \in \R_2 \times \J_1,\label{eq:feas_cut3}
		\end{alignat}\label{prob:M}
	\end{subequations}}
	\noindent where $\J_1$ is the set of all extreme points of
	Problem \eqref{prob:s1} and $\J_2$ and $\R_2$ are the set of all
	extreme points and rays of Problem \eqref{prob:s2},
	respectively.
\end{corollary} 
Let $\C_1$ and $\C_2$ denote the set of all constraints in
\eqref{eq:opt_cut} and \eqref{eq:feas_cut3} respectively.  At each
iteration, the RMP is a relaxation of Problem \eqref{prob:M} with a
subset of the constraints, i.e., $\widetilde \C_1 \subseteq \C_1$
and $\widetilde \C_2 \subseteq \C_2$. The Benders separation routine
at each iteration for an optimal solution $\hat x$ of the RMP is given
by Algorithm \ref{algo:benders} instead of by solving Problem
\eqref{prob:inner:dual} and produces a violated constraints in $\C_i \setminus
\widetilde \C_i,$ for some $i =1,2$.

\begin{algorithm}[!t]
\fontsize{11}{18}\selectfont
\Begin {
	\KwIn{$\hat x \in \mathbb R^{n_1}$}
		Solve Problem \eqref{prob:s1} for $\hat x$;\\
		Obtain its optimal solution $(\hat y, \hat v) \in \J_1$ and let $\Obj$ be its optimal objective value;\\
		Solve Problem \eqref{prob:s2} for $\hat x$ and $\Obj$;\\
			\eIf{Problem \eqref{prob:s2} is unbounded with an unbounded ray $(\tilde \psi, \tilde u_y, \tilde w) \in \R_2$}{
				Add the feasibility cut
				$0 \ge \tilde\psi^T (b- Ax) - \tilde u_y^T (h_y -G_{xy}\hat x) - \tilde w(d^T \hat y  -\hat v^T(k+K_x x))$ to the RMP;}
			{
				Obtain its optimal solution $(\hat \psi, \hat u_y, \hat w) \in \J_2$;\\
			    Add the optimality cut
					$t \ge \hat\psi^T (b- Ax) - \hat u_y^T (h_y -G_{xy}\hat x) - \hat w(d^T \hat y -\hat v^T(k+K_x x))$ to the RMP;
				
					Update the best primal bound with the obtained feasible solution;
			}
		
	}
	\caption{The Benders Separation Algorithm.}
\label{algo:benders}
\end{algorithm}

\subsection{Interpretation of Benders Cuts}\label{sec:interp}
While the Benders cuts (i.e., Equations
\eqref{eq:opt_cut}-\eqref{eq:feas_cut3}) are valid for any $\hat x$
feasible to the SOCP relaxation of Problem (MISOCP), they allow for an
intuitive interpretation when $\hat x \in \X$ (i.e., when the
integrality condition is met). Recall that Remark \ref{rem:theo1:feas} 
indicates, for $\hat x \in \X$,
\begin{subequations}
	\begin{alignat}{2}
	\Obj &= \min_{y \in \K_y, s \ge 0} \{d^Ty + \bar\psi^T s: By + s \ge b - A\hat x\},\label{eq:interp:mibp:s1}\\
	\Obj_{\eqref{prob:s2}} &=\min_{y \in \K_y} \{c_y^T y:By \ge b - A \hat x, \
G_y^T y \ge h_y - G_{xy} \hat x,\ d^T y \le \Obj\}.\label{eq:interp:mibp:s2}
	\end{alignat}\label{eq:interp:mibp}
\end{subequations}
\noindent
Note that Equations \eqref{eq:interp:mibp} imply that, for $\hat x \in
\X$, Problem \eqref{prob:s1} corresponds to the follower problem with a penalty term,
while Problem \eqref{prob:s2} represents the leader's problem
conditional on the follower's reaction, since any $y \ge 0$ satisfying the first and third constraints of Problem \eqref{eq:interp:mibp:s2} is optimal to the follower. 

As upper bounds $\bar \psi$ on the dual follower variables are not available in most cases, each entry of $\bar \psi$ is usually set as a sufficiently large numerical value. Therefore, when the follower problem is infeasible for $\hat x$, some entry of $s$ must take some positive value, incurring a significant cost $\bar \psi^T s$ in the objective function. Note that in that case, Problem \eqref{eq:interp:mibp:s2} becomes infeasible due to the first constraint, generating a cut in the form of \eqref{eq:feas_cut3}, separating the \emph{bilevel-infeasible} $\hat x$. When the follower problem has a finite optimum $\Obj$ for given $\hat x$ (i.e., $s=0$), the first and the third constraints in Problem \eqref{eq:interp:mibp:s2} gurantee that $y$ is feasible to Problem \eqref{eq:interp:mibp:s2} only when $y \in \mathcal F(\hat x)$, i.e., an optimal solution to the follower for given $\hat x$. Therefore, the unboundedness of Problem \eqref{prob:s2} (i.e., Problem \eqref{eq:interp:mibp:s2} is infeasible) implies that none of the follower's response to $\hat x$ is feasible to the leader and the cut \eqref{eq:feas_cut3} correctly cuts off the \emph{bilevel-infeasible} point $\hat x$. Lastly, if both Problems \eqref{prob:s1} and \eqref{prob:s2} have a finite optimum at $\hat x
\in \X$, this implies that $\hat x$ is \emph{bilevel-feasible} and the optimality cut (i.e., Equation \eqref{eq:opt_cut})
correctly evaluates the leader's cost incurred by the follower reaction $\hat y$ (i.e., $f(\hat x)=c_y^T \hat y$).

\section{Numerically Stable Benders Cut Generation Procedure}\label{sec:stable}
Note that, in many cases, there may not be a specific upper bound $\bar \psi$ on $\psi$ available to use in Equation \eqref{eq:Mc}. In those cases, we may use sufficiently large numerical values to set up $\overline \psi$. However, these large values of $\bar \psi$ are undesirable as they may lead to a significantly wide range of coefficients, which not only negatively affects the computation of Problem \eqref{prob:s1} but also compromises the effectiveness of the Benders cuts of Equations \eqref{eq:opt_cut} and \eqref{eq:feas_cut3}. This situation may get worse as the degree of coupling between the leader and the follower problems gets higher.

To address this issue, we propose a numerically stable Benders separation procedure that combines no-good and Benders cuts. Let $g_{\hat x}(x)$ be the function that satisfies $g_{\hat x}(x) = 0$ for $x=\hat x$ and $g_{\hat x}(x)>0$ for any $x \neq \hat x$. If the domain of $g$ is $\mathbb B^{n_x}$, we can define such $g_{\hat x}$ as: $g_{\hat x}(x)=\sum_{j \in \J:\hat x_j = 1} (1-x_j) + \sum_{j\in \J:\hat x_j = 0} x_j$. The idea is to use $g_{\hat x} (x)$ to eliminate the term $\hat v^T (k + K_x x)$ in the Benders cuts, which is associated with the McCormick relaxation (i.e., the term involving an upper bound on $\psi$); it is motivated by a structural role of the term $\hat v^T (k + K_x x)$ that assigns a cost associated $\hat \psi$ to $x$ that deviates from $\hat x$. 

Suppose the RMP generates a guess $(\hat x, \hat t)$. Our objective is to generate a valid inequality that cuts off $(\hat x, \hat t)$ if $\hat x$ is bilevel-infeasible or if $\hat x$ is bilevel-feasible but $(\hat x, \hat t)$ violates Constraint \eqref{prob:bl:1st:constr}. Refer to Definition \ref{def:bilevel_feasible} for the definition of bilevel-feasible points used in this paper. 

The modified procedure first solves the follower problem (i.e., Problem \eqref{prob:bl:lower}) with given $\hat x$. If the follower problem is infeasible, we can cut off the point by using the dual unbounded ray $\tilde \psi_1$ of the follower at $\hat x$: 
\begin{equation}\tilde \psi_1^T(b-Ax) \le 0.
	\label{eq:feas_cut:follower-feasibility}
\end{equation}
Note that this inequality is valid, since any $x$ that violates this inequality is infeasible to the follower, hence bilevel-infeasible. If the follower has an optimal solution $\hat y$ with a finite optimum $\Obj$ at $\hat x$, the procedure solves Problem \eqref{eq:interp:mibp:s2} with given $\hat x$ and $\Obj$ (i.e., the leader problem for the given follower's reaction). If Problem \eqref{eq:interp:mibp:s2} is infeasible with a dual unbounded ray of ($\tilde \psi_2, \tilde u_y, \tilde w$), it adds 
\begin{equation}
    0 \ge \tilde \psi_2^T (b- Ax) + \tilde u_y^T (h_y -G_{xy} x) - \tilde w\left(d^T \hat y + (M - d^T \hat y) g_{\hat x}(x)\right),
    \label{eq:feas_cut:M}
\end{equation}
where $M$ is an upper bound on the follower objective value which can be obtained by solving Problem \eqref{prob:hpp} with the objective function replaced by $d^T y$, and note that such $M$ is guaranteed to exist under Assumption \ref{assum:lb}. Note that it cuts off the bilevel-infeasible solution $\hat x$, since we have $g_{\hat x}(\hat x) = 0$ and thus $\tilde \psi_2^{T} (b- A\hat x) + \tilde u_y (h_y -G_{xy} \hat x) - \tilde w \left(d^T \hat y + (M - d^T \hat y) g_{\hat x}(\hat x)\right) = \tilde \psi_2^{T} (b- A\hat x) + \tilde u_y (h_y -G_{xy} \hat x) - \tilde w \Obj > 0$. Also, note that it does not cut off any bilevel feasible solution $x' \neq \hat x$, since for any extreme ray $(\tilde \psi_2, \tilde u_y, \tilde w)$ of the dual of Problem \eqref{eq:interp:mibp:s2}, the following holds:
\begin{align*}
 0 &\ge \tilde \psi_2^{T} (b- Ax') + \tilde u_y^T (h_y -G_{xy} x') - \tilde w \Obj'\\
   & \ge \tilde \psi_2^{T} (b- Ax') + \tilde u_y^T (h_y -G_{xy} x') - \tilde w M,\\
   & \ge \tilde \psi_2^{T} (b- Ax') + \tilde u_y^T (h_y -G_{xy} x') - \tilde w \left(d^T \hat y + (M - d^T \hat y) g_{\hat x}(x')\right),
\end{align*}
where $\Obj'$ denote the optimal objective of the follower for the given bilevel-feasible $x'$.

Similarly, if Problem \eqref{eq:interp:mibp:s2} is feasible with optimal objective value $\Obj_2 > \hat t$, it adds the following cut:
\begin{equation}
    t \ge \hat \psi_2^T (b- Ax) + \hat u_y^T (h_y -G_{xy} x) - \hat w\left(d^T \hat y + (M - d^T \hat y) g_{\hat x}(x)\right),
    \label{eq:opt_cut:M}
\end{equation}
where ($\hat \psi_2, \hat u_y, \hat w$) is the dual solution of Problem \eqref{eq:interp:mibp:s2}. Note that Equation \eqref{eq:opt_cut:M} cuts off $(\hat x, \hat t)$; since at $\hat x$, we have $g_{\hat x}(\hat x) = 0$ and thus
\[\hat \psi_2 ^T (b- A \hat x) + \hat u_y^T (h_y -G_{xy} \hat x) - \hat w(d^T \hat y) = \Obj_2 > \hat t.\]
For any other bilevel-feasible $x' \neq \hat x$, let $\Obj'$ be the optimal objective value of Problem \eqref{eq:interp:mibp:s1} for given $x'$ and ($\hat \psi_2', \hat u_y', \hat w'$) be the optimal dual solution of Problem \eqref{eq:interp:mibp:s2} for given $x'$ and $\Obj'$. Then, we have
{\fontsize{10}{17}\selectfont
\begin{align*}
 t \ge \hat \psi_2^{'T} (b- Ax') + \hat u_y^{'T} (h_y -G_{xy} x') - \hat w' \Obj'  &\ge \hat \psi_2^{T} (b- Ax') + \hat u_y^T (h_y -G_{xy} x') - \hat w \Obj' \\
    &\ge \hat \psi_2^{T} (b- Ax') + \hat u_y^T (h_y -G_{xy} x') - \hat w M \\
    & \ge \hat \psi_2^{T} (b- Ax') + \hat u_y^T (h_y -G_{xy} x') - \hat w \left(d^T \hat y + (M - d^T \hat y) g_{\hat x}(x')\right)
\end{align*}}
Therefore, it does not cut-off any valid bilevel-feasible $x'$. The modified algorithm is summarized in Algorithm \ref{algo:benders_nogood}.

Note that the benefit of the modified procedure is huge when we do not have a specific upper bound on the dual variables, which happens in many practical problems; in the cut generation procedure, we do not have to deal with large coefficients that may be needed to account for unbounded dual variables. Furthermore, $M$ can be dynamically reduced inside a callback function depending on the relaxation problem at the current branching node as the algorithm proceeds.

\begin{algorithm}[!t]
\fontsize{11}{18}\selectfont
\Begin {
	\KwIn{$\hat x \in \X$}
		Solve Problem \eqref{prob:bl:lower} for $\hat x$;\\
		\eIf{Problem \eqref{prob:bl:lower} is infeasible with a dual unbounded ray $\tilde \psi_1$}{
		    Add the feasibility cut 
				$\tilde \psi_1^T(b-Ax) \le 0$
				to the RMP;
		}{Obtain its optimal solution $\hat y$ and let $\Obj$ be its optimal objective value;\\
		 Solve Problem \eqref{eq:interp:mibp:s2} with $\Obj$ and $\hat x$;\\
			\eIf{Problem \eqref{eq:interp:mibp:s2} is infeasible with a dual unbounded ray $(\tilde \psi_2, \tilde u_y, \tilde w)$}{
				Add the feasibility cut
				$0 \ge \tilde\psi_2^T (b- Ax) - \tilde u_y^T (h_y -G_{xy}\hat x) - \tilde w(d^T \hat y + (M - d^T \hat y) g_{\hat x}(x))$ to the RMP;}
			{
				Obtain its dual optimal solution $(\hat \psi_2, \hat u_y, \hat w)$;\\
			    Add the optimality cut
					$t \ge \hat\psi_2^T (b- Ax) - \hat u_y^T (h_y -G_{xy}\hat x) - \hat w(d^T \hat y + (M - d^T \hat y) g_{\hat x}(x))$ to the RMP;
				
					Update the best primal bound with the obtained feasible solution;
			}
		}
	}
	\caption{The Numerically Stable Benders Separation Algorithm.}
\label{algo:benders_nogood}
\end{algorithm}

\subsection{A Relatively-Complete Follower}\label{sec:relatively-complete-follower}
Instead of ensuring the follower feasibility of $x$ using the feasibility cut \eqref{eq:feas_cut:follower-feasibility}, we may enforce the follower feasibility in the leader problem, which is often shown to be effective in the context of stochastic programming. To make the master problem generates $x$ that is feasible to the follower and likely to be feasible to the leader, we add the follower variable $y$ to the master problem, along with the leader and the follower constraints. Then, we add an additional constraint $t \ge c_y^Ty$ to the master. Note that these additions do not alter the optimal solution and make the initial master problem become equivalent to the HPP problem. The resultant initial master problem is:
\begin{equation}
	\min_{x \in \X, y \in \mathcal K_y} \left\{ c_x^T x + t : t \ge c_y^T y,\ G_{xy} x + G_y y \ge h_y, \
	 Ax + By \ge b\right\},
	\label{prob:initial-master}
\end{equation}
which generates follower-feasible incumbent solutions and eliminates the need to add the feasibility cut \eqref{eq:feas_cut:follower-feasibility}. We will use this extended master problem for the remainder of this paper.

\section{A Heuristic For Finding MIBSOCP Incumbent Solutions}\label{sec:heuristic}
In this section, we propose a heuristic method for finding an incumbent solution of MIBSOCP, which utilizes callback functions available in many commercial branch-and-cut solvers. Callback functions enable users to alter the solvers behavior, such as maneuvering branching rules/directions, adding (lazy) constraints on an as-needed basis while the solver is in process, and updating an incumbent solution.

For example, a user can write a callback function that provides the solver with Benders cuts in a lazy manner only when a newly-found incumbent solution violates some of the Benders cuts. To be specific, a branch-and-cut solver begins to solve the initial master problem without any Benders cuts; and then whenever an incumbent solution $\hat x$ is found in the solver process, it invokes a user-written callback function that checks whether there is a Benders cut violated by $\hat x$; if there exists, it adds the cut to cut off $\hat x$.

Users also can control the branching rules/decisions. Unless a user defines a callback function that alters the branching decision---which node to branch on and how to branch the selected node---the solver, by default, often chooses an integer variable to branch on. For bilevel problems, however, it may not be effective since the follower suboptimality can be a more critical factor in the solution infeasibility than the solution non-integrality.
The following motivating example illustrates the case:
\begin{example}
	Consider the following BLP instance:
	\begin{align}
		\min \ & x - 8y\\
		\mbox{s.t.} \ & 0 \le x \le 11 \mbox{ integer},\\
		& y \in \arg \min_{y \ge 0} \left\{ y : 3x + 4y \ge 18, \ -4x + 9y \le 19,\ 8x+y \le 88\right\},
	\end{align}
	the bilevel-feasible region of which is illustrated in Figure \ref{fig:branching-rule} (a) as a set of filled circles. 

	\begin{figure}[!t]
		\centering
		\subfloat[Graphical representation of the HPP problem: (i) a shaded region: initial master problem feasible region (ii) thick lines: bilevel-feasible region (iii) filled circles: BSOCP feasible solutions (iv) an unfilled circle: relaxation solution at the root node ]{\includegraphics[width=0.5\textwidth]{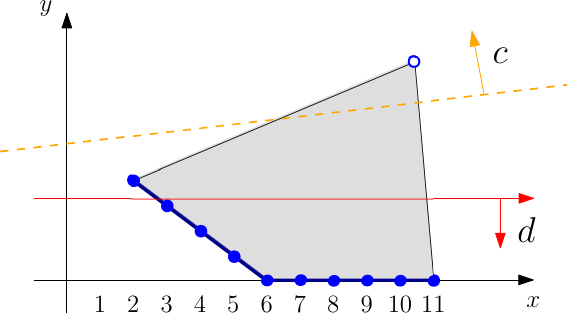}}\\
		  \subfloat[Branching on fractional $x$]{\includegraphics[width=0.4\textwidth]{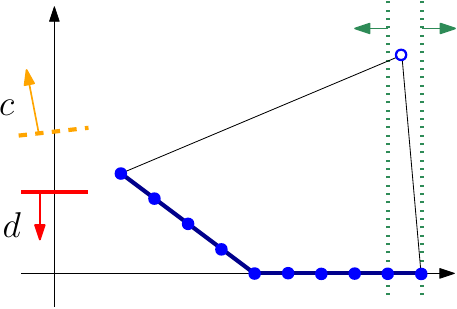}}\hspace*{0.1\textwidth}
		  \subfloat[Proposed branching scheme]{\includegraphics[width=0.4\textwidth]{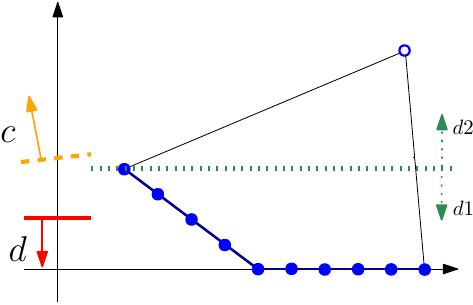}}
		\caption{Motivating Example and Branching Rules for HPP\label{fig:branching-rule}}
		\end{figure}
		\label{exam:branching}
\end{example}

At the root node of a branching tree, the initial master problem (i.e., the HPP problem) consists of the leader and the follower constraints, excluding the integrality condition, which corresponds to the shaded region in Figure \ref{fig:branching-rule} (a). The solution at the root node is $(\frac{773}{76},\frac{126}{19})$ denoted by an unfilled circle in Figure \ref{fig:branching-rule} (a). Figure \ref{fig:branching-rule} (b) illustrates how a solver would branch, by default, on a fractional solution. Note that, for this example, this way of branching neither notably enhances the lower bound nor finds a bilevel-feasible solution. 
This motivates us to explore a bilevel-tailored branching scheme. 

\subsection{A Bilevel-Tailored Branching Rule}\label{sec:heuristic:branching}

A branching rule may be used to guide the solver to a more relevant solution space based on the follower's reaction in the optimization process. Suppose a branch-and-cut solver begins to solve the initial master problem (i.e., Problem \eqref{prob:initial-master}) and adds lazy Benders cuts progressively. Let $(\hat x, \hat y, \hat t)$ be a fractional solution obtained at a branching node. Let $y'$ be the follower's reaction to $\hat x$. If $c_y^T \hat y < c_y^T y'$, then $(\hat x, \hat y)$ is not bilevel-feasible; it is mainly due to the discrepancy between $c_y$ and $d$. The leader's optimal follower response $\hat y$ obtained based on $c_y$ may significantly differ from the actual follower response, i.e., $d^T \hat y > d^T y'$. 

Therefore, a branching scheme on the follower variables that informs the solver of the follower's response may help discover a relevant solution space faster. The idea is that we branch using hyperplanes (d1) $d^T y \le \lfloor\frac{d^Ty' + d^T \hat y}{2} \rfloor$ and (d2) $d^T y > \lfloor\frac{d^Ty' + d^T \hat y}{2} \rfloor$.  We can use this branching scheme when the lower bound does not improve notably after some predetermined number of integer branchings.

This is illustrated in Figure \ref{fig:branching-rule} (c). In the example, unlike the conventional branching rule, a good (essentially optimal) incumbent solution will be found in the node generated by (d1). Note also that, although the branching direction (d2) does not have any bilevel feasible solutions, it will keep branching on in the direction. This unnecessary exploration can be prevented if we can feed the solver a proper follower upper bound. For example, the follower objective value in Example 1 cannot be greater than 3, and thus if the algorithm knows $d^Ty \le 3$, it will prune the direction (d2) immediately. Then, the question is how we can efficiently obtain a meaningful upper bound of the follower problem, which motivates Section \ref{sec:heuristic:local-cuts}.

\subsection{Local Cuts}\label{sec:heuristic:local-cuts}
A meaningful upper bound of the follower may be obtained in various ways. Let FUB denote an upper bound of the follower; then, a cut $d^T y \le \mbox{FUB}$ may inform the solver of a bilevel-infeasible region. Note that FUB may decrease as we move from the root node to a leaf node of a search tree, and thus making the cut stronger. Therefore, if we can find FUB at each branching node efficiently, a local cut $d^T y \le \mbox{FUB}$ (i.e., a cut that is valid for a node to which it is added and all nodes descending from the node) may prune a direction with no follower feasible solution like (d2) in Example \ref{exam:branching}.

\cite{fischetti2017new} proposed a convenient way of finding a follower upper bound that utilizes the branching decisions progressively made in branch-and-bound algorithms. FUB is obtained by solving a restricted follower problem in which the leader variables are fixed at values that most restrict the follower problem in each branching node. However, this scheme may not give a meaningful upper bound as the restriction is often too strict. 

On the other hand, the least upper bound of the follower at each branching node can be obtained by solving the following max-min problem:
\begin{subequations}
\begin{alignat}{3}
	\max_{x_j \in [\underline x_j', \overline x_j']_{\mathbb Z}, \forall j \in \mathcal J, x \in HPP_x} \ & \min_{y \in \mathcal K_y} \left\{ d^Ty: By \ge b-Ax\right\}, 
\end{alignat}\label{prob:follower-ub}
\end{subequations}
where $\underline x_j'$ and $\overline x_j'$ respectively represents the lower and the upper bounds on $x_j$ at the current branching node and $HPP_x$ denotes the projection of the feasible region of the HPP problem onto the space of $x$. Note that Problem \eqref{prob:follower-ub} can be solved by replacing the inner minimization problem with a set of linear inequalities containing the reversed weak duality constraint and the McCormick relaxation as illustrated in Section \ref{sec:MIP}. However, the resultant problem is a MISOCP problem, which may be hard to solve efficiently. 

Several approximations/relaxations of Problem \eqref{prob:follower-ub} may be used to improve its solution time. In this paper, we use an approximation of FUB at each branching node by finding a fixed point $x'$ (i.e., an equilibrium point). Consider a relaxation solution $\hat x$ at the current branching node. First, we solve the inner minimization problem with $x$ fixed at $\hat x$, the result of which affects the objective value of the outer problem by ${\hat \psi}^T (b-A\hat x)$, where $\hat \psi$ is the dual solution of the inner problem. Based on this updated objective information, the outer problem may respond with a different solution $\hat x'$; we obtain the outer problem's next action by solving $\max_{x_j \in [\underline x_j', \overline x_j'], \forall j \in \mathcal J, x \in \overline{HPP}_x}-{\hat \psi}^TAx$, where $HPP_x$ is replaced with its LP relaxation $\overline{HPP}_x$ to improve the computation time. This procedure is repeated until we find a fixed point $x'$ (i.e., $x'$ remains unchanged after an iteration); we use the objective value at $x'$ as FUB.

As an approximation may produce an invalid FUB, a local cut obtained by this restriction may not be valid. Therefore, we use the scheme proposed in this section as a heuristic method for finding a good incumbent solution along with an initial set of Benders cuts in the preprocessing step, which is illustrated in Algorithm \ref{algo:benders:preprocessing}. It is worth noting, however, that if a valid relaxation of Problem \eqref{prob:follower-ub} is used, the method can be used in the overall process of the solver.

\begin{algorithm}[!t]
	\fontsize{11}{18}\selectfont
	\Begin {
		Turn on the user-written callback functions for the bilevel-tailored branching scheme and the heuristic local cut in the numerically stable benders method;\\
		Solve the instance within some predetermined time limit (e.g., 150 sec);	
		
		Turn off the bilevel-tailored branching scheme and the heuristic local cut generation;

			\eIf{Heuristic cut has been added}{
			Resolve the instance from scratch with the incumbent solution and Benders cuts found in the preprocessing steps;}{
				Resume the solver process;
			}
		}
		\caption{The Numerically Stable Benders Separation Algorithm With A Preprocessing Step.}
	\label{algo:benders:preprocessing}
	\end{algorithm}

\section{BSOCP with Additional Upper Level Constraints on Dual Variables of Lower Level}\label{sec:ext}

An interesting extension of BSOCP is to add an additional constraint
to the upper level problem which states the impact of the follower
dual variables on the leader problem. Such constraints may be
desirable in some sequential market environment where the follower
dual variables settle the prices of commodities that are used by the
leader, e.g., a unit-commitment problem for interdependent
natural gas and electricity markets studied by
\cite{byeon2019unit}. This section discusses how the BSP decomposition
technique carries over to this extension.

In order to formulate the situation where the follower's dual solution affects the leader problem, BSOCP can be extended as follows:
\begin{subequations}
\begin{alignat}{4}
 \quad &\min_{x,y,\psi} \qquad && c_x^T x + c_y^T y\label{eq:bl2:obj}\\
& \mbox{ s.t.}    && G_{xy} x + G_y y \ge h_y, \label{eq:bl2:upper:y}\\
&   && G_{x\psi} x + G_\psi \psi \ge h_\psi, \label{eq:bl2:upper:psi}\\
& &&  x \in \X:=\{x \in \K_x: x_i \in \mathbb B, \ \forall i \in \I\},\label{eq:bl2:upper:x}\\
& && (y, \psi) \in \Q \left(\min_{y \in \K_y} \{d^T y: A x + B y \ge b \}\right), \label{prob:bl2:lower}
\end{alignat}\label{prob:bl2}
\end{subequations}
\noindent
where $\psi$ denote dual variables of the follower, $\Q(P)$ denotes the set of optimal primal and dual solution pairs of Problem $P$, and $G_{x\psi}, G_\psi, h_\psi$ are given rational matrices or vectors of appropriate dimension.

For this extension, the definition of bilevel-feasibility can be naturally extended as follows:
\begin{definition} 
	A leader decision $\hat x$ is called \emph{bilevel-feasible}, if it satisfies all of the following: 
	\begin{enumerate}
		\item $\hat x$ is feasible to the follower (i.e., Problem \eqref{prob:bl:lower} with $x$ fixed as $\hat x$ is feasible).
		\item there exists an optimal response $(\hat y, \hat \psi)$ of the follower for given $\hat x$ that is feasible to the leader, i.e., $G_{xy} \hat x + G_y \hat y \ge h_y$, $G_{x\psi} \hat x + G_y \hat y \ge h_\psi$, and $(\hat y, \hat \psi) \in \mathcal F(\hat x)$, where $\mathcal F(\hat x)$ denote the set of optimal primal and dual solution pairs of the follower problem at $\hat x$.
	\end{enumerate}
	When at least one of the above is not met, we call the pair \emph{bilevel-infeasible}.
	\label{def:bilevel_feasible:2}
\end{definition}

Let Problem (MISOCP)$'$ denote the MISOCP reformulation of Problem \eqref{prob:bl2}. It is easy to see that Problem (MISOCP)$'$ is equivalent to Problem (MISOCP) to which Constraint \eqref{eq:bl2:upper:psi} added. Let $u_\psi$ denote the dual variable associated with Constraint \eqref{eq:bl2:upper:psi}. Then, the dual of Problem \eqref{prob:inner} with Constraint \eqref{eq:bl2:upper:psi} (i.e., the BSP for (MISOCP)$'$) is expressed as Problem \eqref{prob:inner:dual} with additional terms $u_\psi^T (h_y - G_{x\psi} \hat x)$ on the objective and $-G_\psi^T u_\psi$ on the left-hand side of Constraint \eqref{eq:inner:dual:primal1}, which we call Problem \eqref{prob:inner:dual}$'$.
\begin{theorem}
	The BSP for (MISOCP)$'$ can be solved by solving two more tractable problems sequentially, i.e., solve the problem
	\begin{subequations}
	\begin{alignat}{2}
	\min_{y \in \K_y,u_\psi \ge 0,v \ge 0} \ &  d^T y - u_\psi^T (h_\psi - G_{x\psi} \hat x) - v^T (k+K_x\hat x)  \\
    \mbox{ s.t.}  \ & By -  G_\psi^T u_\psi - K_\psi^T v \ge b, \label{eq:s1:2:1}\\
     & K_\mu^T v \le a,\label{eq:s1:2:2}
	\end{alignat}\label{prob:s1:2}
	\end{subequations}
	and then solve Problem \eqref{prob:s2} where $$\Obj=\begin{cases}
		\mbox{Optimal objective value of Problem \eqref{prob:s1:2} } & \mbox{if Problem \eqref{prob:s1:2} has a finite optimum}\\
		\infty & otherwise.
	\end{cases}$$
\label{theo:bender:2}
\end{theorem}

\begin{remark}
Note that Problem \eqref{prob:s1:2} is a relaxation of Problem \eqref{prob:s1} with an additional vector of variables $u_\psi$, so Problem \eqref{prob:s1:2} is also guaranteed to be feasible by Remark \ref{rem:theo1:feas}. Consider the dual of Problem \eqref{prob:s1:2}: 
\[\max_{\psi\ge 0,\mu \ge 0}\left\{b^T\psi - \mu^T a : B^T \psi \preceq_{\K_y} d, \ G_\psi \psi \ge h_\psi - G_{x\psi}\hat x,\ K_\psi \psi + K_\mu \mu \ge k + K_x \hat x\right\},\]
	where $\psi$ and $\mu$ are dual variables associated with Constraints \eqref{eq:s1:2:1} and \eqref{eq:s1:2:2} respectively. For any $\hat x \in \X$, the McCormick relaxation is exact and the optimal objective value of the dual of Problem \eqref{prob:s1:2} becomes equivalent to 
	\begin{equation}
		\max_{0 \le \psi \le \bar \psi} \{\psi^T (b - A\hat x) : B^T \psi \preceq_{\K_y} d, G_\psi \psi \ge h_\psi - G_{x\psi}\hat x\}. \label{prob:s1:2:dual:hatx}
	\end{equation} 
	Note that Problem \eqref{prob:s1:2:dual:hatx} can be infeasible, as it can fail to satisfy $G_\psi \psi \ge h_\psi - G_{x\psi}\hat x$. Therefore Problem \eqref{prob:s1:2} can be unbounded for some $\hat x$.
	
	\label{rem:theo2}
\end{remark}

As a result of Theorem \ref{theo:bender:2}, Corollary \ref{coro:1} also extends to this case as follows:
\begin{corollary}\label{coro:1:2}
	Problem (MISOCP)$'$ is equivalent
	to the following problem:
	{\fontsize{8}{17}\selectfont
		\begin{subequations}
			\begin{alignat}{3}
			\min_{x \in \X} \ & c_x^T x + t && \nonumber\\
				\mbox{ s.t.} \ & t \ge \hat\psi^T (b- Ax) + \hat u_y^T (h_y -G_{xy} x) - \hat w\left(d^T \hat y - \hat u_\psi^T(h_\psi - G_{x\psi} x) -\hat v^T(k+K_x x)\right), &&\forall (\hat \psi, \hat u_y, \hat w,\hat y, \hat u_\psi, \hat v) \in \J_2 \times \J_1, \label{eq:dual:opt_cut}\\
			& d^T \tilde y - \tilde u_\psi^T(h_\psi - G_{x\psi} x) -\tilde v^T(k+K_x x) \ge 0, &&\forall (\tilde y, \tilde u_\psi, \tilde v) \in \R_1, \label{eq:dual:feas_cut}\\
				& 0 \ge \tilde\psi^T (b- Ax) + \tilde u_y^T (h_y -G_{xy} x) - \tilde w\left(d^T \hat y - \hat u_\psi^T(h_\psi - G_{x\psi} x) -\hat v^T(k+K_x x)\right), \ && \forall (\tilde \psi, \tilde u_y, \tilde w,\hat y,\hat u_\psi,\hat v) \in \R_2 \times \J_1,\label{eq:dual:feas_cut3}
			\end{alignat}\label{prob:M:dual}
		\end{subequations}}
		\noindent where $\J_1$ and $\R_1$ are the set of all extreme points and rays of
		Problem \eqref{prob:s1:2} and $\J_2$ and $\R_2$ are the set of all
		extreme points and rays of Problem \eqref{prob:s2},
		respectively. 
\end{corollary}

The Benders separation routine for $\hat x$ is given by Algorithm \ref{algo:benders:2}.

\begin{algorithm}[!t]
	\fontsize{11}{18}\selectfont
	\Begin {
		\KwIn{$\hat x \in \mathbb R^{n_x}$}
			Solve Problem \eqref{prob:s1:2};\\
			\eIf{Problem \eqref{prob:s1:2} is unbounded with an unbounded ray $(\tilde y, \tilde u_\psi, \tilde v)$}{
				Add the feasibility cut 
					$d^T \tilde y - \tilde u_\psi^T (h_\psi - G_{x\psi}x) -\tilde v^T(k+K_x x ) \ge 0$
					to the RMP;
			}{Obtain its optimal solution $(\hat y, \hat u_\psi, \hat v)$ and let $\Obj$ be its optimal objective value;\\
			Solve Problem \eqref{prob:s2} for $\hat x$ and $\Obj$;\\
			\eIf{
				Problem \eqref{prob:s2} is unbounded with an unbounded ray $(\tilde \psi, \tilde u_y, \tilde w)$
			}{
				Add the feasibility cut
					$0 \ge \tilde\psi^T (b- Ax) + \tilde u_y^T (h_y -G_{xy} x) - \tilde w\left(d^T \hat y - \hat u_\psi^T(h_\psi - G_{x\psi} x) -\hat v^T(k+K_x x)\right)$
			}{
				
				Obtain its optimal solution $(\hat \psi, \hat u_y, \hat w)$;\\
				Add the optimality cut
					$t \ge \hat\psi^T (b- Ax) + \hat u_y^T (h_y -G_{xy} x) - \hat w\left(d^T \hat y - \hat u_\psi^T(h_\psi - G_{x\psi} x) -\hat v^T(k+K_x x)\right)$ to the RMP;
				
				Update the best primal bound with the obtained feasible solution;
			}	
		}
	}
	\caption{The Benders Separation Algorithm for (MISOCP)$'$.}
	\label{algo:benders:2}
\end{algorithm}

\subsection{Interpretation of Benders Cuts}\label{sec:interp_dual}
Equations
\eqref{eq:dual:opt_cut}-\eqref{eq:dual:feas_cut3} also allow for an intuitive interpretation of the Benders cuts for $\hat x \in \X$. Recall that Remark \ref{rem:theo2} 
indicates, for $\hat x \in \X$, the optimal objective value of Problem \eqref{prob:s1:2} equals that of Problem \eqref{prob:s1:2:dual:hatx}. Let $\Obj_{\eqref{prob:bl:lower}}$ denote the optimal objective value of the follower at $\hat x$. 

First, consider the case where Problem \eqref{prob:s1:2} is unbounded for given $\hat x \in \X$, which implies the infeasibility of Problem \eqref{prob:s1:2:dual:hatx} for $\hat x$. Note that this means that there is no follower dual solution that satisfy $G_\psi \hat\psi \ge h_\psi - G_{x\psi}\hat x$. Therefore, it adds the cut \eqref{eq:dual:feas_cut} to cut off the \emph{bilevel-infeasible} point $\hat x$.

When Problem \eqref{prob:s1:2} has a finite optimum $\Obj$ for given $\hat x \in \X$, there are two possible scenarios: (i) $\Obj = \Obj_{\eqref{prob:bl:lower}}$ or (ii) $\Obj < \Obj_{\eqref{prob:bl:lower}}$. Note that when Case (ii) is the case, Problem \eqref{eq:interp:mibp:s2} must be infeasible (i.e., Problem \eqref{prob:s2} is unbounded), since $\Obj_{\eqref{prob:bl:lower}}$ is the smallest possible value of $d^Ty$ while satisfying $By \ge b-A\hat x$. Note that Case (ii) occurs when there is no follower \emph{optimal} dual solution $\hat \psi$ at $\hat x$ that satisfies the constraint (i.e., $G_\psi \hat\psi < h_\psi - G_{x\psi}\hat x, \forall (\hat y, \hat \psi) \in \mathcal Q(\hat x)$). Any dual optimal solution $\hat \psi$ is no longer feasible to Problem \eqref{prob:s1:2:dual:hatx} and needs to be altered so that $G_\psi \psi \ge h_\psi - G_{x\psi}\hat x$ becomes satisfied, lowering the optimal objective value of Problem \eqref{prob:s1:2:dual:hatx} than that of the follower problem at $\hat x$. Therefore, Case (ii) implies that the leader constraint on the follower's dual response is not met. Therefore, $\hat x$ is \emph{bilevel-infeasible}, so the cut \eqref{eq:dual:feas_cut3} cuts off $\hat x$.

When Case (i) is the case (which implies that there exists the follower's dual response $\hat \psi$ to $\hat x$ that satisfies the leader constraint on the follower's dual response), the dual of Problem \eqref{prob:s1:2:dual:hatx} can be considered as the same as the follower problem \eqref{prob:bl:lower}, hence the same interpretation as in Section \ref{sec:interp} holds.  

\begin{remark}
	Note that, for this general case, it may be difficult to derive the numerically stable benders cut generation procedure. In order to obtain numerically stable cuts that are valid for any bilevel feasible solution, a reasonable value of $M$---an upper bound on Problem \eqref{prob:s1:2}---is needed, as in Equations \eqref{eq:feas_cut:M} and \eqref{eq:opt_cut:M}. However, as discussed in Remark \ref{rem:theo2}, Problem \eqref{prob:s1:2} can be unbounded, and thus for this general case we cannot rely on the numerically stable Benders decomposition. 
	\label{rem:general:numerically_stable}
\end{remark}

Accordingly, we propose a special class of BSOCP that allows for a sequence-independent decomposition as well as several acceleration schemes that can be applied to improve the computational performance of Algorithm \ref{algo:benders:2} in Sections \ref{sec:sequence-independent} and \ref{sec:accel}.

\subsection{Sequence-Independent BSP Decomposition}
\label{sec:sequence-independent}

Some special cases of BSOCP allow for a stronger alternative to
Theorem \ref{theo:bender:2}. In this section, we deal with the extended
version of BSOCP discussed in Section \ref{sec:ext} (i.e., Problem
\eqref{prob:bl2}), but any result in this section also holds for
Problem \eqref{prob:bl}. As noted in Section \ref{sec:ext}, the BSP
of Problem (MISOCP)$'$ is decomposed into two problems, i.e., Problems
\eqref{prob:s1:2} and \eqref{prob:s2}, which are solved in a
sequential manner. A sequence-independent BSP decomposition is allowed
in two special cases of BSOCP: (i) $d=c_y$ (ii) $c_y = 0$. Case (i)
subsumes a class of mixed-integer conic-linear optimization problems
that involves additional constraints on the dual variables of its
inner-continuous problem, which is the case of \cite{byeon2019unit}.

\begin{corollary}
Let Problem \eqref{prob:s2}$'$ denote Problem \eqref{prob:s2} with $w$ fixed at zero. Then, the BSP for Problem (MISOCP)$'$ with $d = c_y$ can be solved by solving
  Problem \eqref{prob:s1:2} and Problem \eqref{prob:s2}$'$ \textbf{independently}.
  \label{coro:independent1}
\end{corollary}
\noindent

A similar result holds for Case (ii).
\begin{corollary}
Let Problem \eqref{prob:s2}$''$ denote Problem \eqref{prob:s2} with $w$ fixed at zero and the right-hand side of Equation \eqref{eq:s2:1} replaced with $d$. Then, the BSP for Problem (MISOCP)$'$ with $c_y = 0$ can be solved by solving
  Problem \eqref{prob:s1:2} and Problem \eqref{prob:s2}$''$ \textbf{independently}.
  \label{coro:independent2}
\end{corollary}

Corollary \ref{coro:independent1} (or \ref{coro:independent2}) implies that the Benders cuts for BSOCP with $d = c_y$ (or $c_y = 0$) can be obtained by solving Problems \eqref{prob:s1} and \eqref{prob:s2}$'$ (or \eqref{prob:s2}$''$) independently and comparing their objective values; This simplifies the Benders cut generation algorithm as described in Algorithm \ref{algo:independent}.

\begin{algorithm}[!t]
	\fontsize{11}{15}\selectfont
\Begin {
	\KwIn{$\hat x \in \mathbb R^{n_1}$}
		Solve Problems \eqref{prob:s1:2} and \eqref{prob:s2}$'$ (or \eqref{prob:s2}$''$) independently and let $\Obj_1$ and $\Obj_2$ respectively denote their objective value;\\
		\eIf{$\Obj_1 = -\infty$ with an unbounded ray $(\tilde y, \tilde u_\psi, \tilde v)$}{
				Add the feasibility cut 
				$c_y^T \tilde y - \tilde u_{\psi}^T(h_\psi - G_{x\psi} x) - \tilde v^T(k+K_x x) \ge 0$
				to the RMP;}
		{
			\eIf{$\Obj_2 = \infty$ with an unbounded ray $(\tilde \psi, \tilde u_y)$}{

				Add the feasibility cut 
				$\tilde \psi^T (b- Ax) + \tilde u_y^T (h_y - G_{xy}\hat x) \le  0$
				to the RMP;
				
			}{
				Obtain the optimal solution $(\hat y, \hat u_\psi, \hat v)$ of Problem \eqref{prob:s1:2};\\
				\eIf{$\Obj_1 < \Obj_2$}{
					Add the feasibility cut 
					\begin{equation}\hat \psi^T (b- Ax) + \hat u_y^T (h_y - G_{xy}\hat x) \le c_y^T \hat y - \hat u_{\psi}^T(h_\psi - G_{x\psi} x) - \hat v^T(k+K_x x)\label{eq:feas_cut:seq}\end{equation}
					to the RMP;}
				{
					Add the optimality cut
					\begin{equation}t \ge \hat \psi^T (b- Ax)+ \hat u_y^T (h_y - G_{xy}\hat x) \mbox{ (or $t \ge 0$)}\label{eq:opt_cut:seq} \end{equation} to the RMP;
					Update the best primal bound with the obtained feasible solution;
				}
			}
		}
	}
	\caption{The Benders Separation Method for BSOCP with $d=c_y$ (or $c_y=0$).}
\label{algo:independent}
\end{algorithm}

\begin{remark}
	The Benders cuts obtained using Algorithm \ref{algo:independent} also allow for an intuitive interpretation. Note that for both Case (i) and (ii), Problem \eqref{prob:s2}$'$ and Problem \eqref{prob:s2}$''$ become equivalent to the following problem:
	$$\Obj_{\eqref{prob:s2}} = \min_{y \in \K_y} \{d^Ty: By \ge b-A \hat x, G_y^Ty \ge h_y - G_{xy}\hat x\}.$$
	Also, for $\hat x \in \X$, Problem \eqref{prob:s1:2} is equivalent to the dual of Problem \eqref{prob:s1:2:dual:hatx}:
	$$\Obj = \min_{y \in \K_y, u_\psi \ge 0, s \ge 0} \{d^Ty + u_\psi^T (h_\psi - G_{x\psi}\hat x) + \bar \psi^T s: By -G_{\psi}^Tu_\psi + s\ge b-A \hat x\}.$$

	Note that, by construction, $\Obj \le \Obj_{\eqref{prob:s2}}$ always holds and the equality holds if and only if $\hat x$ is bilevel-feasible. When $\Obj < \Obj_{\eqref{prob:s2}}$, for some $\hat x \in \X$, it implies that either the leader's constraint on the follower's optimal primal solution or that on the follower's dual solution is not satisfied. Therefore, the cut \eqref{eq:dual:feas_cut3} can be replaced by a cut \eqref{eq:feas_cut:seq} that enforces $\Obj \ge \Obj_{\eqref{prob:s2}}$. Also, for any bilevel-feasible $\hat x$ (i.e., with $\Obj \le \Obj_{\eqref{prob:s2}}$), the optimality cut \eqref{eq:opt_cut:seq} correctly evaluates the cost incurred.
	\label{rem:sequence-independent}
\end{remark}

\subsection{Acceleration Schemes}\label{sec:accel}

This section presents some acceleration schemes for the standard Benders decomposition method discussed in previous literature (e.g., \cite{fischetti2010note} and \cite{ben2007acceleration}) and shows that these schemes can be applied to the dedicated Benders method for the general case described in Section \ref{sec:ext}. 

\subsubsection{Normalizing Benders Feasibility Cuts}
\label{sec:norm}

\cite{fischetti2010note} have shown that normalizing
the ray used in Benders feasibility cuts can improve the performance of
Benders decomposition. The Benders subproblem decomposition outlined in Algorithm \ref{algo:benders:2} can be generalized to produce a 
normalized ray.

When Problem \eqref{prob:inner:dual}$'$ is unbounded, the problem at hand consists in solving Problem
\eqref{prob:inner:dual}$'$ to which an additional normalization constraint of $\|(\psi, u_y,w, y, u_\psi, v)\|_1 = 1$ is added and with the right-hand side of other constraints set to zero. Let \texttt{ubd}\eqref{prob:inner:dual}$'$ denote the resultant problem. 
The proof of Theorem \ref{theo:bender:2} showed
that \texttt{ubd}\eqref{prob:inner:dual}$'$ has three different types of extreme rays:
\begin{enumerate}
    \item[(i)] $\tilde \mu_1 := ( 0, 0, 0,\tilde y, \tilde u_\psi,  \tilde v)$ for $(\tilde y, \tilde u_\psi, \tilde v) \in \R_1$.
    \item[(ii)] $\tilde \mu_3 := (\tilde \psi, \tilde u_y, \tilde w,\tilde w\hat y, \tilde w \hat u_\psi, \tilde w\hat v)$ for $(\hat y, \hat u_\psi, \hat v) \in \J_1$ and $(\tilde \psi, \tilde u_y, \tilde w) \in \R_2$ with $\tilde w = 0$.
    \item[(iii)] $\tilde \mu_3 := (\tilde \psi, \tilde u_y, \tilde w,\tilde w\hat y, \tilde w \hat u_\psi, \tilde w\hat v)$ for $(\hat y, \hat u_\psi, \hat v) \in \J_1$ and $(\tilde \psi, \tilde u_y, \tilde w) \in \R_2$ with $\tilde w > 0$.
\end{enumerate}

\noindent
Cases (i) and (ii) are simple: It suffices to solve Problem
\eqref{prob:s1:2} and Problem \eqref{prob:s2} with the additional
constraint of $\|( y, u_\psi, v)\|_1 = 1$ and $\|(\psi, u_y, w)\|_1 = 1$ respectively.
Case of (iii) (i.e., when Problem \eqref{prob:s1:2} has a finite optimum $\Obj$ at $(\hat y,  \hat u_\psi, \hat v) \in \J_1$ and Problem \eqref{prob:s2} is unbounded with an unbounded ray of $(\tilde \psi, \tilde u_y, \tilde w) \in \R_2$ with $\tilde w > 0$ and $\U := \tilde \psi^T (b - A \hat x) + \tilde u_y^T (h_y - G_{xy} \hat x) - \Obj \tilde w > 0$) is more difficult and requires to find a normalized ray
$\tilde r' = (\tilde \psi', \tilde u_y', \tilde w',\tilde y', \tilde u_\psi',\tilde v')$ that
maximizes the objective function of \texttt{ubd}\eqref{prob:inner:dual}$'$ while satisfying $\|\tilde r'\|_1
= 1$ and $\tilde w' > 0$. Note that $\tilde \mu_3 / \|\tilde \mu_3 \|_1$ is a feasible solution
to \texttt{ubd}\eqref{prob:inner:dual}$'$. Hence, \texttt{ubd}\eqref{prob:inner:dual}$'$ is
feasible and bounded. 

Consider the Lagrangian relaxation of \texttt{ubd}\eqref{prob:inner:dual}$'$ with
$w>0$ that penalizes the violation of the normalization constraint
with some $\lambda \in \mathbb R$.
By defining ($\psi, u_y, w, y, u_\psi, v$) = ($\frac{\psi}{w},$ $\frac{u_y}{w}$, $1$, $\frac{y}{w}$, $\frac{u_\psi}{w}$, $\frac{v}{w}$), the Lagrangian relaxation becomes as follows:
\begin{equation}
\min_{\lambda \in \mathbb R}\left\{\lambda + \sup_{w > 0} \left\{ w t(\lambda)\right\}\right\},
\label{prob:LagRelax}
\end{equation}
where $t(\lambda) := t^2(\lambda) - t^1(\lambda) - \lambda$ and 
\begin{multline}
		t^1(\lambda)=  \min_{y\in \K_y,(y^+, y^-, u_\psi, v)^T \ge 0}  \left\{ d^T y + \lambda \textbf 1^T (y^+ + y^-) -u_\psi^T(h_\psi - G_{x\psi}\hat x-\lambda \textbf 1) \right.\\
		- v^T(k+K_x \hat x-\lambda \textbf 1) :
	By  -G_\psi^T u_\psi - K_\psi^T v  \ge b,  \ K_\mu^T v \le a, y= y^+ - y^- \left.\right\},
	\label{prob:t1}
\end{multline}
and 
\begin{equation}
	t^2(\lambda)=\max_{\psi \ge 0, u_y \ge 0} \left\{ \psi^T (b-A\hat x-\lambda \textbf 1) + u_y^T(h_y - G_{xy} \hat x-\lambda \textbf 1):B^T \psi + G_y^T u_y  \le d \right\},\label{prob:t2}
\end{equation}
\begin{proposition}
	An optimal solution $\lambda^*$ of Problem \eqref{prob:LagRelax} is the solution of $t(\lambda) = 0.$\label{prop:normalizing-feasibility-cut}
\end{proposition}
\noindent Since $t(\lambda)$ is a convex piecewise linear function of $\lambda$, the solution of $t(\lambda) = 0$ can be found via a subgradient-based Newton's method as shown in Algorithm
\ref{algo:sub}. At each iteration $k$, $-(\hat \psi^k, \hat u_{y}^k, 1, \hat y^{+k},
\hat y^{-k}, \hat u_{\psi}^k, \hat v^k)^T \textbf 1$, where $(\hat y^{+k},
\hat y^{-k},\hat u_{\psi}^k, \hat v^k)$ and $(\hat \psi^k, \hat u_{y}^k)$ are the solutions of $t^1(\lambda^k)$ and $t^2(\lambda^k)$
respectively, is a subgradient of $t$ at $\lambda^k$ and is denoted by
$\delta t(\lambda^k)$. $\lambda^{k+1}$ is a solution of a linear approximation of $t(\lambda)$ at $\lambda^k$ (i.e., $\lambda^{k+1}$ is the solution of $\delta t(\lambda^k)(\lambda - \lambda^k)+t(\lambda^k)=0)$. Observe that Problems \eqref{prob:t1} and
  \eqref{prob:t2} are the counterparts to Problem \eqref{prob:s1:2} and
  \eqref{prob:s2}, demonstrating that the subproblem decomposition
  carries over to the decomposition.

\begin{algorithm}[!t]
	\fontsize{11}{18}\selectfont
\Begin {
	\KwIn{$\lambda^0 = 0$, $t(\lambda^0) =\frac{\U}{\tilde w}$, $k = 0$;}
	\While{$t(\lambda^k) > \epsilon$}{
	    Calculate $\delta t(\lambda^k)$ (a subgradient of $t$ at $\lambda = \lambda^k$);\\
	    $\lambda^{k+1} = \lambda^k - \frac{t(\lambda^k)}{\delta t(\lambda^k)}$;\\
	    Solve $t^1(\lambda^{k+1})$ and $t^2(\lambda^{k+1})$ and calculate $t(\lambda^{k+1}) = t^2(\lambda^{k+1}) - t^1(\lambda^{k+1}) - \lambda^{k+1}$;\\
	    $k \leftarrow k+1;$\\
	}
	}
\caption{The Subgradient Newton's Method for Problem \eqref{prob:Lag2}.}
\label{algo:sub}
\end{algorithm}


\subsubsection{An In-Out Approach} 
\label{sec:inout}

\cite{ben2007acceleration} proposed an acceleration
scheme (the in-out method) for general cutting-plane algorithms. The
method carefully chooses the separation point, rather than using the
solution obtained from the RMP. The method considers two points: a
feasible point $x_{in}$ to Problem \eqref{prob:M} and the optimal
solution $x_{out}$ of the RMP. It uses a convex combination of these
two points when generating the separating cut, i.e., it solves Problem
\eqref{prob:inner:dual} with $\hat x = \lambda x_{in} + (1-\lambda) x_{out}$
for some $\lambda \in (0,1)$.

\cite{fischetti2016redesigning} applied the in-out
approach with an additional perturbation to solve facility
location problems:
\begin{equation}
    \hat x = \lambda x_{in} + (1-\lambda) x_{out} + \epsilon \textbf 1, \label{eq:sep}
\end{equation}
for some $\lambda \in (0,1)$ and $\epsilon > 0$, and showed a
computational improvement.

This paper also employs the in-out approach equipped with some
perturbation as \cite{fischetti2016redesigning}. It
periodically finds $x_{in}$ in a heuristic manner and chooses the
separation point according to Equation \eqref{eq:sep}. The
implementation starts with $\lambda = 0.5$ and $\epsilon = 10^{-6}$
and decrease $\lambda$ by half if the BD halts (i.e., it does not
improve the optimality gap for more than 3 consecutive iterations). If
the algorithm halts and $\lambda$ is smaller than $10^{-5}$,
$\epsilon$ is set to $0$. After 3 more consecutive iterations without
a lower bound improvement, the algorithm returns to the original BD.
Whenever a new best incumbent solution is found, the in-out approach
is applied again with this new feasible point.

    \section{Computational Results}\label{sec:comp}
This section studies the performance of the numerically stable Benders cut generation procedure, proposed in Section \ref{sec:stable}, along with the performance of the heuristic method in Section \ref{sec:heuristic}, and that of the dedicated Benders method for problems with leader's constraints on follower's dual, proposed in Section \ref{sec:ext}.

\subsection{Performance Analysis of the Numerically Stable Benders Method}\label{sec:result:numerically-stable-algorithm}
\paragraph{Benchmark: \cite{fischetti2017new}.} We compare the proposed algorithm with a publicly available state-of-the-art solver for MIBLP \citep{fischetti2017test}, which implements the algorithm proposed by \cite{fischetti2017new} along with some acceleration schemes. As the benchmark algorithm can only be applied to linear cases, we run experiments only for BLP in this section, and the result on BSOCP will be presented in Section \ref{sec:result:dedicated-benders-algorithm}.

\paragraph{Test cases.} We obtained BLP test instances by relaxing the follower integrality constraints of general bilevel test cases \texttt{XUWANG} proposed by \cite{xu2014exact} and \texttt{XUWANG-LARGE}, \texttt{MIPLIB} produced by \cite{fischetti2017new}; the original instances are available at an open-source repository \citep{fischetti2017test}.

\paragraph{Implementation}: 
The numerically stable Benders method, proposed in Section \ref{sec:stable}, is implemented with the C++/Cplex interface and all the experiments were executed on a virtual Linux machine with 13.5 GB of memory allocated on an Intel Core i7 PC at 2.3 GHz. Each run has a wall-time limit of 1 hour. As the current version of the benchmark algorithm cannot be run on a Cplex with a version higher than 12.7.1, we used the Cplex 12.7.1 library for the benchmark algorithm. For the proposed Benders method, we used Cplex 20.1.0 due to some technical issues in callback functions encountered when implementing the heuristic method with older versions of Cplex. When the heuristic method is not used, we observed that there was no noticeable difference between the results obtained by the proposed method using Cplex 12.7.1 and 20.1.0.

The Benders cuts are implemented using a user-defined callback class, inherited from \texttt{Lazy\-Const\-raint\-Call\-backI}; whenever the master problem finds an incumbent solution $\hat z = (\hat x, \hat y, \hat t)$, the callback class instance is called during the optimization process and checks whether $\hat z$ violates any of the feasibility or optimality cuts by solving the subproblems (Problems \eqref{prob:bl:lower} and \eqref{eq:interp:mibp:s2}); if there exists such a Benders cut that cuts off $\hat z$, it is added to the master problem as a lazy constraint. In addition, if an optimality cut is found to cut off $\hat z$, then $\hat x$ combined with the solution $y'$ to Problem \eqref{eq:interp:mibp:s2} is a bilevel-feasible solution; thus, if its objective value (i.e., $c_x^T\hat x + c_y^T y'$) improves the current upper bound, we use another user-defined callback class, inherited from \texttt{HeuristicCallbackI}, to update the best incumbent solution. 

The heuristic method proposed in Section \ref{sec:heuristic} is also implemented using callback. The bilevel-tailored branching rule discussed in Section \ref{sec:heuristic:branching} is implemented using user-written callback classes, inherited from \texttt{Branch\-Call\-backI} and \texttt{Node\-Call\-backI}. If the best lower bound does not improve for 3 consecutive nodes, it solves the follower problem with given solution $\hat x$ at the current node and obtain its optimal objective value $d^Ty'$. Then, it branches on (i) $d^T y \le \lfloor \frac{d^Ty' + d^T \hat y}{2} \rfloor$ and (ii) $d^Ty > \lfloor \frac{d^Ty' + d^T \hat y}{2} \rfloor$. \texttt{User\-Cut\-Call\-backI} is used to implement the heuristic local cut proposed in Section \ref{sec:heuristic:local-cuts}; at each branching node, the solver calls a user-written function that obtains an approximate FUB and adds a heuristic local cut $d^T y \le$ FUB. For computationally hard instances, the numerical Benders method is initially equipped with the branching and the local cut callbacks for at most \texttt{150} seconds, as in Algorithm \ref{algo:benders:preprocessing}.

Other than callbacks, we set the integrality and feasibility tolerances as \texttt{1e-9}, and other parameters were set as default values. We also applied a preprocessing step proposed in \cite{fischetti2017new} which fixes some of $y$-variables if it is guaranteed to have a fixed value due to the optimality. 

\paragraph{Result on \texttt{XUWANG} and \texttt{XUWANG-LARGE}.}
\begin{figure}[t!]
  \centering
  \includegraphics[width=0.4\textwidth]{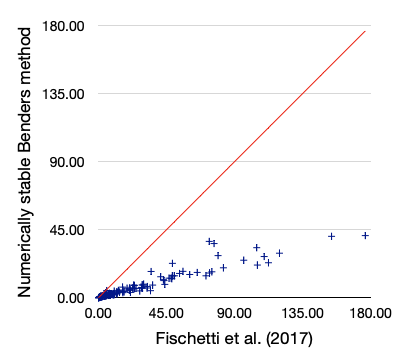}
  \caption{Computational comparison on \texttt{XUWANG} and \texttt{XUWANG-LARGE} instances}\label{fig:result:xuwang}
\end{figure}
The result on \texttt{XUWANG} and \texttt{XUWANG-LARGE} is displayed in Figure \ref{fig:result:xuwang}. The horizontal and vertical axes respectively represent the computation time of \cite{fischetti2017new} and the proposed algorithm without the upper-bounding method in seconds, and each point on the figure shows the computation times taken by the benchmark algorithm (horizontal axis) and the proposed algorithm (vertical axis) for solving an instance. Therefore, points under the red line represent the instances where the proposed algorithm was faster than the benchmark algorithm. Note that, for all instances, the proposed algorithm solved the instances about 3 times faster than the benchmark algorithm on average.

\paragraph{Result on \texttt{MIPLIB}.} Table \ref{table:result:miplib} in Appendix \ref{appendix:result:miplib} summarizes the computational performance of \cite{fischetti2017new}, the numerically stable Benders method (denoted by \texttt{N}), and \texttt{N} equipped with the heuristic proposed in Section \ref{sec:heuristic} (denoted by \texttt{N+H}). Note that the heuristic method improves the optimality gap of \texttt{N} in many instances, notably for air03-5, cap6000-1, cap6000-5, cap6000-9, harp2-1, and harp2-5. 

For \texttt{N+H} and the benchmark algorithm, the instances that \texttt{N+H} outperforms are denoted by boldface letters, i.e., better optimality gap or better computation time when their optimality gaps tie. It is hard to tell which one is better; the benchmark algorithm performed notably better on air04-9, air05-9, cap6000-1, enigma-5, harp2-9, and \texttt{N+H} does on cap6000-5, cap6000-9, harp2-5, nw4-5, nw4-9, p0201-1, and p0201-5. It is worth noting that the Benders cut itself does not take account of integrality of the master problem, thus the Benders cuts may not be as tight as the cuts that employ the integrality information of the problem. This may be the reason why the Benders method takes more time to close the gap for enigma-5 than the benchmark algorithm does, which avails intersection cuts, and did not solve cap6000-1 to the optimality while the benchmark algorithm does, even though \texttt{N+H} finds the optimal solution in an early stage. On the other hand, \texttt{N+H} gives incumbent solutions on p0201-1 and p0201-5 for which the benchmark algorithm suffers from numerical issues and finds a better incumbent solution on the instances that it outperforms. The result suggests that a hybrid of the proposed method and the benchmark algorithm may improve the computation time further.

\subsection{Performance Analysis of the Dedicated Benders Method proposed in Section \ref{sec:ext}}\label{sec:result:dedicated-benders-algorithm}
This section studies the performance of the decomposition approach
(Section \ref{sec:decom}) and the benefits of the acceleration schemes
explained in Sections \ref{sec:norm} and \ref{sec:inout}.  All
algorithms were implemented with the C++/Gurobi 8.0.1 interface and
executed on an Intel Core i5 PC at 2.7 GHz with 8 GB of RAM. Each run
has a wall-time limit of 1 hour.

\subsubsection{Test Instances} 
\label{sec:test}
A recent paper by \cite{byeon2019unit}
introduced the unit commitment problem with Gas Network Awareness
(UCGNA), a tri-level optimization problem where the first and second
levels determine how to commit and dispatch electric power generating
units; The third level decides how to operate the gas network given
the natural gas demands of committed gas-fueled generators that are
determined in the first and second levels. The economic feedback from
the gas network, i.e., the natural gas zonal prices, is given by the
dual solution $\psi$ of the third-level optimization and the first-level
optimization is subject to constraints over both $\psi$ and commitment
decisions $x$ in order to ensure the robustness of the unit commitment
decisions against the economic feedback from the gas system. \cite{byeon2019unit} showed that the tri-level problem can be
reformulated as a special case of BSOCP discussed in Section \ref{sec:sequence-independent}. The detail of the model is given in Appendix \ref{sec:UCGNA}.
The evaluation of the proposed method is performed on the instances of the UCGNA problem.

The instances are based on the gas-grid test
system, which is representative of the natural gas and electric power
systems in the Northeastern United States \citep{bent2018joint}. There
are 42 different instances, each of which constructed by uniformly
increasing the demand of each system by some percentage; $\eta_p$
denotes the stress level imposed on the power system which takes
values from $\{1, 1.3, 1.6\}$ and $\eta_g$ denotes the stress level of
the gas system that has values of $\{1,1.1,\cdots, 2.2,2.3\}$. For
example, $(\eta_p, \eta_g) = (1.3, 2.3)$ means the demands of the
power and natural gas systems are increased uniformly by $30\%$ and
$130\%$ respectively. Before we experiment with the solution approaches on the instances of the UCGNA problem, we apply some preprocessing step which eliminates invalid bids with regard to a lower bound on natural gas zonal prices. A detailed description of the instances and the preprocessing step can be
found in \citep{byeon2019unit}.

\subsubsection{Computational Performance}
\label{sec:comp:three}
This section compares three different solution approaches for BSOCP:
\begin{itemize}
\item[\texttt{D:}] the proposed dedicated Benders method with the acceleration schemes (Section \ref{sec:ext});
\item[\texttt{G:}] an off-the-shelve solver (Gurobi 8.0.1);
\item[\texttt{B:}] the standard Benders method with the acceleration
  schemes (Section \ref{sec:accel}).
\end{itemize}
 
The implementation of $\texttt{D}$ is sequential, although Problems
\eqref{prob:s2}$'$ and \eqref{prob:s1:2} can be solved independently (See Corollary \ref{coro:independent1}). All
solution approaches use the same values for the Gurobi parameters,
i.e., the default values except \texttt{NumericFocus} set at 3,
\texttt{DualReductions} at 0, \texttt{ScaleFlag} at 0,
\texttt{BarQCPConvTol} at 1e-7, and \texttt{Aggregate} at 0 for more
rigorous attempts to detect and manage numerical issues.
 
Tables \ref{table:computation1}-\ref{table:computation16} in Appendix \ref{appendix:result:UCGNA} report the
computation times and optimality gaps of the three solution
methods. The symbol $\dagger$ indicates that a method reaches the time
limit and the symbol $\ddagger$ that the method did not find any
incumbent solution. The results for $\eta_p = 1$ are summarized in Table
\ref{table:computation1}; \texttt{D} timed out for two instances,
\texttt{G} reached the time limit for 5 instances, and \texttt{B}
timed out for all the instances. For the two instances with $\eta_g =
1.8, 1.9$, where all methods time out, \texttt{D} found incumbent
solutions within optimality gaps of 1.8\% and 1.3\% and \texttt{B}
found solutions with gaps of 6.7\% and 10.6\%. On the other hand,
\texttt{G} did not find any incumbent solution. For easy instances
that both \texttt{D} and \texttt{G} found optimal solutions within two
minutes, \texttt{G} is faster than \texttt{D} by a factor of 2 in
average.

For instances with $\eta_p = 1.3$, reported in Table
\ref{table:computation13}, \texttt{D} and \texttt{G} timed out for 7
instances and \texttt{B} reached the time limit for all the
instances. For the 7 instances with $\eta_g = 1.6, \cdots, 2.2$, where
all methods reached the time limit, \texttt{D} found incumbent
solutions within 4.3\% of optimality and \texttt{B} found worse
solutions. On the other hand, \texttt{G} did not find any incumbent
solution except the two instances with $\eta_g = 1.6$ and $2$. For
easy instances that both \texttt{D} and \texttt{G} found optimal
solutions within two minutes, \texttt{G} is faster than \texttt{D} by
a factor of around 7 in average.

Instances with $\eta_p = 1.6$ display similar behaviors. While
\texttt{B} failed to find optimal solutions for all the instances,
\texttt{D} and \texttt{G} found optimal solutions for 7 instances. For
the hard instances where all methods timed out, \texttt{D} found
incumbent solutions with optimality gaps less than 7.5\%, \texttt{B}
found worse solutions, and \texttt{G} failed to find any incumbent
solution. For the instances where both \texttt{D} and \texttt{G} found
optimal solutions, \texttt{G} is faster than \texttt{D}.

\begin{figure}[!t]
\centering
\subfloat[Computation Time (sec).]{\includegraphics[width=0.4\textwidth]{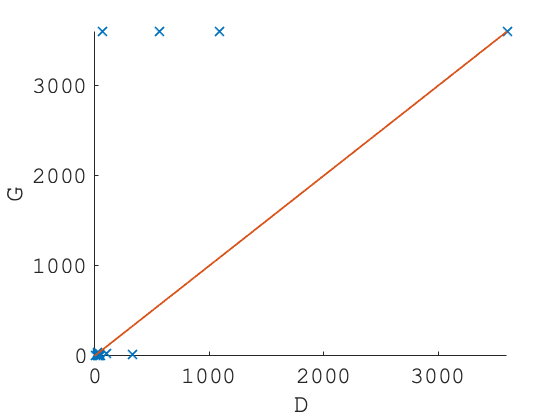}\label{fig:DvsG:time}}
\subfloat[Optimality Gap (\%, logarithmic scale).]{\includegraphics[width=0.4\textwidth]{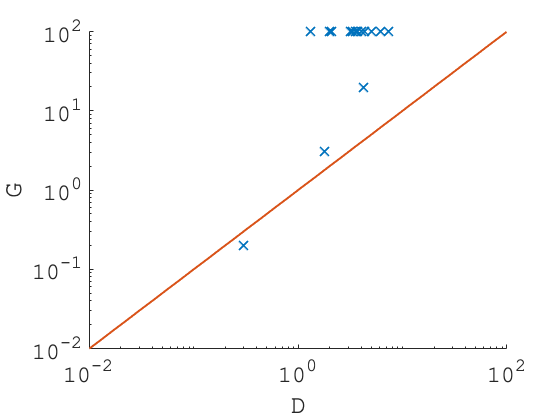}\label{fig:DvsG:gap}}
\caption{\texttt{D} vs \texttt{G}.\label{fig:DvsG}}
\end{figure}

To compare the computational performance of \texttt{D} and \texttt{G}
more precisely, Figure \ref{fig:DvsG} visualizes the performance of
\texttt{D} and \texttt{G} for all the instances. Figure
\ref{fig:DvsG:time} reports the computation times of \texttt{D} and
\texttt{G}, Figure \ref{fig:DvsG:gap} displays the optimality gaps of
the two methods for all the instances, and the reference lines (in
red) serve to delineate when a method is faster than the other. For
Figure \ref{fig:DvsG:gap}, the axes are in logarithmic scale and a
100\% optimality gap is assigned to instances with no incumbent. The
figure indicates that, although \texttt{D} is slower than \texttt{G}
for some easy instances (the points at the bottom left corner of
Figure \ref{fig:DvsG:time}), it has notable benefits for hard
instances (the points in the upper left side of Figures
\ref{fig:DvsG:time} and \ref{fig:DvsG:gap}).

\subsubsection{Benefits of the Acceleration Schemes}

This section studies the benefits of the acceleration schemes by
comparing the performance of the dedicated Benders method with
different combinations of acceleration schemes applied. It uses
\texttt{D($n_k$,$i_k$)} to denote the dedicated Benders method 
with acceleration schemes ($n_k$,$i_k$) where
\begin{itemize}
\item $n_k$: $k = 1$ if the normalization scheme is applied; $k = 0$ otherwise; 
\item $i_k$: $k = 1$ if the in-out approach is applied; $k = 0$ otherwise.
\end{itemize}

Tables \ref{table:accel1}-\ref{table:accel16} in Appendix \ref{appendix:result:UCGNA} summarize the
computational performance of the dedicated Benders methods with the
four combinations of acceleration schemes. Table \ref{table:accel1} displays the computation times and optimality
gaps for instances with $\eta_p = 1$. Without the in-out approach,
\texttt{D}$(n_1, i_0)$ and \texttt{D}$(n_0, i_0)$ timed out for all
instances. Although both \texttt{D}$(n_1, i_0)$ and \texttt{D}$(n_0,
i_0)$ reach the time limit for all instances, the normalization scheme
does improve optimality gaps. On the other hand, with the in-out
approach, \texttt{D}($n_0, i_1$), solves 10 instances within 100
seconds. However, \texttt{D}($n_0, i_1$) still cannot solve the two
instances with $\eta_g = 2.1, 2.2$. The slight increase in computation
time of \texttt{D}($n_1, i_1$) for some instances, compared to
\texttt{D}($n_0, i_1$), is due to the additional computation time
required to find a normalized ray.

The results for instances with $\eta_p = 1.3$ are reported in Table
\ref{table:accel13}. Again, without the in-out approach,
\texttt{D}$(n_1, i_0)$ and \texttt{D}$(n_0, i_0)$ timed out for all
instances, but \texttt{D}$(n_1, i_0)$ has significant improvement in
optimality gaps for some instances. With the in-out approach,
\texttt{D}$(n_0, i_1)$ solved 7 instances within 150 seconds and so
did \texttt{D}$(n_1, i_1)$. The normalization scheme does have some
computational benefits, as \texttt{D}$(n_1, i_1)$ has smaller
optimality gaps than \texttt{D}$(n_0, i_1)$ for the remaining 7
instances except one instance with $\eta_g = 2.2$. Moreover, for some
hard instances where \texttt{D}$(n_0, i_1)$ reached the time limit,
\texttt{D}$(n_1, i_0)$ has smaller optimality gaps (i.e., $\eta_g =
1.7, \cdots, 2$).

The acceleration schemes display similar behaviors for instances with
$\eta_p = 1.6$. Without the in-out approach, \texttt{D}$(n_0, i_0)$
timed out for all instances, while \texttt{D}$(n_1, i_0)$ solves
one instance to optimality and has significant improvements in
optimality gaps. With the in-out approach, both \texttt{D}$(n_0, i_1)$
and \texttt{D}$(n_1, i_1)$ solve 7 instances within 350 seconds, and
\texttt{D}$(n_1, i_1)$ has smaller optimality gaps for the unsolved
instances. Again, for some hard instances for which \texttt{D}$(n_0, i_1)$
reached the time limit, \texttt{D}$(n_1, i_0)$ has smaller optimality
gaps (i.e., $\eta_g = 1.7, \cdots, 2.2$).

\subsubsection{Benefits of the Decomposition Method}

Section \ref{sec:comp:three} indicated that the decomposition method
has significant benefits for solving BSOCP. The decomposition method
not only shortens computation times required for solving the dual
of the inner-continuous problem, but also allows us to address the
numerical issues of BSOCP.

\begin{figure}[!t]
\centering
\includegraphics[width=0.45\textwidth]{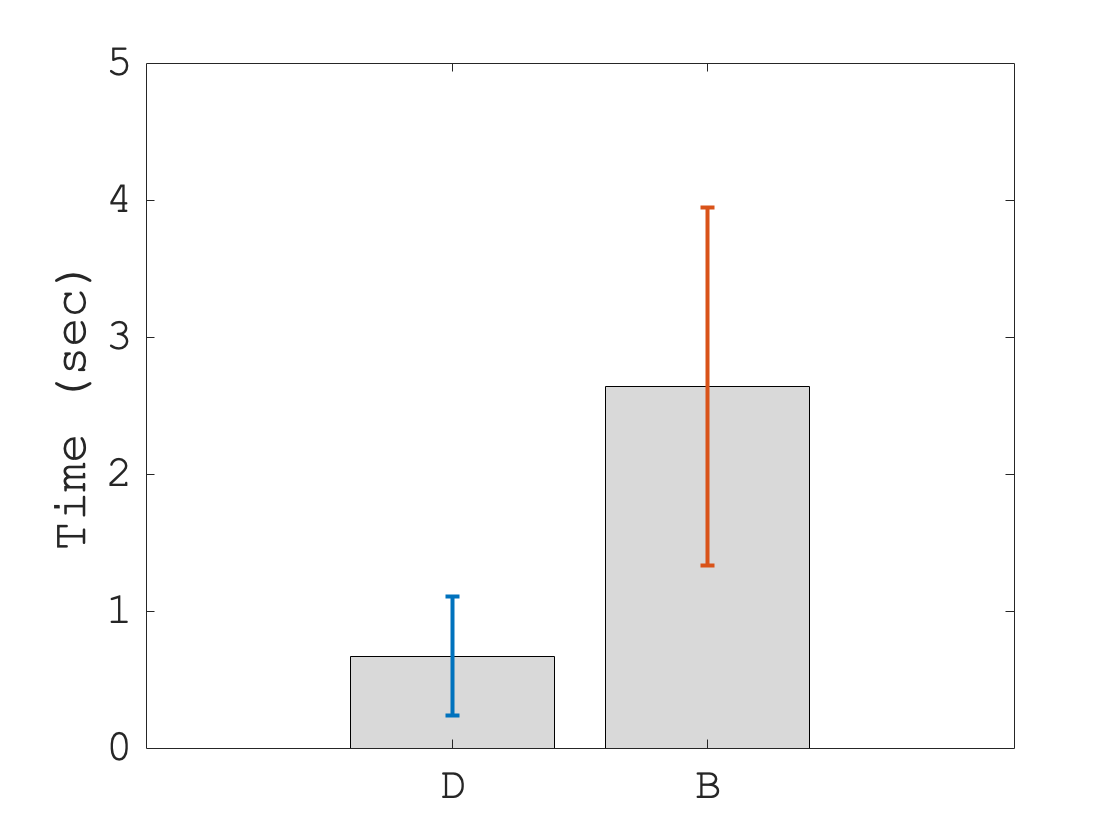}
\caption{Statistics on Computation Times for Cut Generation.\label{fig:iter_time}}
\end{figure}

Figure \ref{fig:iter_time} displays the average computation time for
generating a Benders cut, where the error bars represent the standard
deviation. On average, the cut generation time of \texttt D is faster
than \texttt B by a factor of 3.94. Since the subproblems that \texttt
D solves to generate cuts (i.e., Problems \eqref{prob:s2}$'$ and
\eqref{prob:s1}) can be solved independently, implementation in
parallel computing would improve the computation time even further.

Moreover, the decomposition method deals better with numerical issues
arising from the complex inner-continuous problem of BSOCP. Figure
\ref{fig:conv} in Appendix \ref{appendix:result:UCGNA} displays the convergence behavior of \texttt{D} and
\texttt{B} for two instances, ($\eta_p, \eta_g) = (1, 1.2), (1.6,
1.8)$. For instance ($\eta_p, \eta_g) = (1, 1.2)$ (i.e., Figure
\ref{fig:conv1D} and Figure \ref{fig:conv1B}), \texttt D closes the
gap in 30 seconds, but \texttt B does not improve its lower bound even
if it finds a good incumbent solution early. For instance ($\eta_p,
\eta_g) = (1.6, 1.8)$ (i.e., Figure \ref{fig:conv2D} and Figure
\ref{fig:conv2B}), although both \texttt D and \texttt B timed out,
\texttt B improves its lower bound much slower than \texttt D. This
behavior of \texttt B is explained by the fact that it suffers from
numerical issues when solving Problem \eqref{prob:inner:dual}; it sometimes
terminates with an optimal solution even if there exists an unbounded
ray. This incorrect evaluation of the first-stage variable leads to
ineffective cut generation and a slower convergence rate. On the other
hand, the decomposition method effectively decomposes Problem
\eqref{prob:inner:dual} into two more stable and smaller problems, which
addresses the numerical issues effectively.

	
    \section{Conclusion}
\label{sec:conclusion}

BSOCP is an important class of hierarchical optimization models that
arises in many practical contexts, including network planning/design
problems in energy systems and transportation networks, facility
location problems, and unit scheduling problems under interdependent
markets. This paper proposed a dedicated Benders decomposition algorithm to
solve BSOCP models, recognizing that the Benders subproblem cannot
necessarily be solved efficiently for large BSOCP problems. The dedicated
approach decomposes the Benders subproblem into two more tractable,
sequentially solvable problems that are closely related to the leader
and the follower problems. The paper showed that the Benders subproblem
decomposition can also be applied to the extension of BSOCP where the
upper-level problem has additional constraints on the leader variables
and the follower dual variables. The paper also discussed a couple of subclasses
of BSOCP that allows the subproblems to be solved independently. In addition, the paper showed how to (i) generate numerically stable cuts by hybridizing no-good and Benders cuts, (ii) obtain a good incumbent solution, and (iii) adapt existing acceleration schemes to this decomposition. In particular, the paper  (i) eliminated the need of arbitrarily large big-M values arising from McCormick reformulations, (ii) made use of novel branching decisions and local cuts that exploit the nature of bilevel optimization, and lastly (iii) showed how to normalize Benders feasibility cuts using a Newton's (subgradient) method and how to carefully choose the separation points using the in-out approach \citep{ben2007acceleration}.

The proposed method was compared with a state-of-the-art bilevel-tailored branch-and-cut algorithm \cite{fischetti2017new}, demonstrating the benefits of the numerically-stable cuts and the heuristic method on BLP instances. The result suggests a potential benefit of a hybrid use of Benders and intersection cuts for bilevel problems, which is left to future research. Also, the proposed decomposition significantly improves the performance of
a standard Benders method and outperforms a state-of-the-art
mathematical-programming solvers for hard BSOCP instances. The experimental
results highlighted the benefits of acceleration schemes---normalizing
feasibility rays and the in-out approach---and demonstrated that
decomposing the Benders subproblem not only shortens the computation
time for generating Benders cuts but also addresses the numerical
issues arising when solving complex Benders subproblems.

\section*{Acknowledgment}

This research was partly supported by an NSF CRISP Award (NSF-1638331)

	\bibliographystyle{informs2014.bst} 
	\bibliography{references.bib}
	
	
	\newpage
	\begin{APPENDICES}
\section{Proof of Theorem \ref{theo:bender}.}
	\label{sec:proof:theo}
	The proof strategy is to
	show that there is a surjective mapping from the possible outcomes of
	Problems \eqref{prob:s1} and \eqref{prob:s2} to those of Problem
	\eqref{prob:inner:dual}, which implies that Problem \eqref{prob:inner:dual} is
	completely determined by Problems \eqref{prob:s1} and \eqref{prob:s2}.
	
	Let $U_{(i)}$ and $F_{(i)}$ respectively denote the unbounded and finite outcome
	of Problem ($i$) for $i \in \{\ref{prob:inner:dual}, \ref{prob:s1},
	\ref{prob:s2}\}$. Due to Remark \ref{rem:theo1:feas}, the combination of all possible outcomes of
	Problems \eqref{prob:s1} and \eqref{prob:s2} are given by 
	\[
	\A = \left\{ 
		(F_{\eqref{prob:s1}},U_{\eqref{prob:s2}}),
		(F_{\eqref{prob:s1}},F_{\eqref{prob:s2}})
		\right\}.
	\]
	Likewise, the possible outcomes of Problem \eqref{prob:inner:dual} can
	be expressed as $\B=\{ U_{\eqref{prob:inner:dual}}, F_{\eqref{prob:inner:dual}}\}$. The
	proof gives a surjective mapping $g:\A \rightarrow \B$, showing the
	solution of Problem \eqref{prob:inner:dual} can be obtained from the solutions
	of Problems \eqref{prob:s1} and \eqref{prob:s2}.

	Let $(\hat y, \hat v)$
	be the optimal solution of Problem \eqref{prob:s1} and $\Obj$ denotes its optimal objective value. 
	
	
	
		\begin{enumerate}
			\item 
			Outcome $U_{\eqref{prob:s2}}$: 
			Let $(\tilde \psi, \tilde u_y, \tilde w)$ denote the unbounded ray of Problem \eqref{prob:s2}. Note that $\tilde \mu_3 := (\tilde \psi, \tilde u_y, \tilde w,\tilde w \hat y, \tilde w \hat v)$ is a feasible ray to Problem \eqref{prob:inner:dual} and has a positive objective value of $\U_{\eqref{prob:s2}}:=\tilde \psi^T (b - A \hat x) + \tilde u_y^T (h_y - G_{xy} \hat x) - \Obj w > 0.$ Therefore $\tilde \mu_3$ is an unbounded ray of Problem \eqref{prob:inner:dual} and Problem \eqref{prob:inner:dual} is unbounded.
		
		\item Outcome $F_{\eqref{prob:s2}}$: Let $(\hat \psi, \hat u_y, \hat w)$ denote the optimal solution of Problem \eqref{prob:s2} and denote its optimal objective value as $\Obj_{\eqref{prob:s2}}$. 
		
		The proof is by a case analysis over two
		versions of Problem \eqref{prob:inner:dual} in which $w > 0$ and $w =
		0$. Note first that $\hat \mu := (\hat \psi, \hat u_y, \hat w,\hat w \hat y,\hat w \hat v)$ is a feasible solution to Problem \eqref{prob:inner:dual} and has an objective value of $\Obj_{\eqref{prob:s2}}$. Suppose $w>0$, then by stating $(y, v) = (\frac{y}{w},  \frac{v}{w})$,
		Problem \eqref{prob:inner:dual} becomes as follows:
		\begin{equation}
			\max_{w > 0} \Obj(w),
			\label{prob:inner:dual:w}
		\end{equation}
		where
		\begin{subequations}
		\begin{alignat*}{4}
			\Obj(w) :=  &\max \ && \psi^T (b - A \hat x) + u_y^T (h_y - G_{xy} \hat x) -w\Obj\\
		& \mbox{ s.t.}  &&  B^T \psi + G_y^T u_y  \preceq_{\K_y} dw + c_y,\\
		& &&\psi \ge 0,u_y \ge 0.
		\end{alignat*}
		\end{subequations}
		
		Note that Problem \eqref{prob:inner:dual:w} is equivalent to Problem \eqref{prob:s2} where the nonnegativity constraint for $w$ is restricted by strict inequality. Therefore, $\max_{w > 0} \Obj(w) \le \Obj_{\eqref{prob:s2}}.$
		
		When $w = 0$, Problem \eqref{prob:inner:dual} can be decomposed into Problem \eqref{prob:s2} with $w$ fixed at 0 (i.e., a restriction of Problem \eqref{prob:s2}) and
		\begin{equation}
		\min \{d^T y - v^T (k+K_x\hat x) : By  - K_\psi^T v \ge 0,\ K_\mu^T v \le 0,\ y \in \K_y,v \ge 0\}.
			\label{prob:s10}
		\end{equation}
		Note that Problem \eqref{prob:s10} is either unbounded or zero at optimality,
		since it has a trivial solution with all variables at
		zeros. Therefore, its optimum must be zero since otherwise Problem
		\eqref{prob:s1} is unbounded. This implies that the optimal objective value of Problem \eqref{prob:inner:dual} when $w=0$ is also bounded above by $\Obj_\eqref{prob:s2}$, which proves that $\hat \mu$ is the optimal solution of Problem \eqref{prob:inner:dual}.
		\qed
	
	\end{enumerate}

\section{Proof of Corollary \ref{coro:1}.}
The proof of Theorem \ref{theo:bender} implies that $\hat \mu$ is an extreme point of Problem \eqref{prob:inner:dual} if and only if $\hat \mu = (\hat \psi, \hat u_y, \hat w,\hat w \hat y,\hat w \hat v)$ for some $(\hat \psi, \hat u_y, \hat w,\hat y, \hat v) \in \J_2 \times \J_1$. Therefore, Equation \eqref{eq:opt_cut} holds. Likewise, the proof of Theorem \ref{theo:bender} also indicates that $\tilde \mu$ is an extreme ray of Problem \eqref{prob:inner:dual} if and only if $\tilde \mu=(\tilde \psi, \tilde u_y, \tilde w,\tilde w\hat y, \tilde w\hat v)$ for $(\hat y, \hat v) \in \J_1$ and $(\tilde \psi, \tilde u_y, \tilde w) \in \R_2$. Thus, Equation \eqref{eq:feas_cut3} holds. This implies that Equations \eqref{eq:opt_cut} and \eqref{eq:feas_cut3} are equivalent to Constraint \eqref{prob:bl:1st:constr}.  \qed

\section{Proof of Theorem \ref{theo:bender:2}.}
	\label{sec:proof:theo2}
	The proof strategy is similar to Theorem \ref{theo:bender:2}: to
	show that there is a surjective mapping from the possible outcomes of
	Problems \eqref{prob:s1:2} and \eqref{prob:s2} to those of the BSP of (MISOCP)$'$, denoted by Problem (BSP).
	
	Let $U_{(i)}$ and $F_{(i)}$ respectively denote the unbounded and finite outcome
	of Problem ($i$) for $i \in \{\mbox{BSP}, \ref{prob:s1:2},
	\ref{prob:s2}\}$. Due to Remark \ref{rem:theo2}, the combination of all possible outcomes of
	Problems \eqref{prob:s1:2} and \eqref{prob:s2} are given by 
	\[
	\A = \left\{ U_{\eqref{prob:s1:2}},
		(F_{\eqref{prob:s1:2}},U_{\eqref{prob:s2}}),
		(F_{\eqref{prob:s1:2}},F_{\eqref{prob:s2}})
		\right\}.
	\]
	Likewise, the possible outcomes of Problem (BSP) can
	be expressed as $\B=\{ U_{(BSP)}, F_{(BSP)}\}$. The
	proof gives a surjective mapping $g:\A \rightarrow \B$, showing the
	solution of Problem (BSP) can be obtained from the solutions
	of Problems \eqref{prob:s1:2} and \eqref{prob:s2}.
	
	\begin{enumerate}
	\item Outcome $U_{\eqref{prob:s1:2}}$: Let $(\tilde y, \tilde u_\psi, \tilde v)$ be the unbounded ray of Problem \eqref{prob:s1:2} and $\U := d^T\tilde y -\tilde u_\psi^T(h_\psi - G_{x\psi} \hat x) - \tilde v (k + K_x \hat x) < 0$. Note that, by construction, $\Obj = \infty$, and thus we can assume w.l.o.g., $w=0$. Note that the feasibility of Problem \eqref{prob:s2} with $w = 0$ is guaranteed due to Assumption \ref{assum:lb}. Let $(\psi', u_y', 0)$ and $\Obj' < \infty$ respectively be any feasible solution of Problem \eqref{prob:s2} and its corresponding objective value. 
	Then, for any $\alpha >0$, $(\psi', u_y', 0, 0, 0, 0) + \alpha ( 0, 0, 0,\tilde y, \tilde u_\psi, \tilde v)$ is a feasible solution to Problem (BSP) and has an objective value of $\Obj' - \alpha \U,$
	which increases as $\alpha$ increases. Hence $\tilde \mu_1=( 0, 0, 0,\tilde y, \tilde u_\psi, \tilde v)$ is an unbounded ray of Problem (BSP) and Problem (BSP) is unbounded.
	
	\item Outcome $F_{\eqref{prob:s1:2}}$: Similar arguments as in Theorem \ref{theo:bender} apply.
	\qed

	\end{enumerate}	

\section{Proof of Corollary \ref{coro:independent1}.} 
Built upon Theorem \ref{theo:bender:2}, it suffices to show that solving Problem \eqref{prob:s2}$'$ is sufficient to obtain the optimal solution or unbounded ray of Problem \eqref{prob:s2}. Note that, by defining $(\psi',u_y') = (\frac{\psi}{w+1},\frac{u_y}{w+1})$, Problem \eqref{prob:s2} becomes as follows:
	\begin{subequations}
	\begin{alignat}{2}
	\max \ & {\psi'}^T (b - A \hat x) + {u_y'}^T (h_y - G_{xy} \hat x) + w \left[{\psi'}^T (b - A \hat x)+ {u_y'}^T (h_y - G_{xy} \hat x)-\Obj\right]  \\
    \mbox{ s.t.}  \ &B^T \psi' + G_y^T u'_y \le c_y, \label{eq:independent:s2:1}\\
     & \psi',u'_y, w \ge 0.
	\end{alignat}\label{prob:independent:s2}
	\end{subequations}
Suppose Problem \eqref{prob:s2}$'$ has a finite optimum $\Obj_{\eqref{prob:s2}'}$ at $(\hat \psi, \hat u_y)$ and $\Obj_{\eqref{prob:s2}'} > \Obj$. Then, for any $\alpha >0$, $(\psi', u_y', w) = (\hat \psi, \hat u_y, \alpha)$ is feasible to Problem \eqref{prob:independent:s2} and its objective value increases as $\alpha$ increases, and thus Problem \eqref{prob:independent:s2} is unbounded, so is Problem \eqref{prob:s2}. Note that by converting $(\psi', u_y', w)$ to the solution of Problem \eqref{prob:s2} using $(\psi',u_y') = (\frac{\psi}{w+1},\frac{u_y}{w+1})$, we can see that $(\hat \psi, \hat u_y, 1)$ is an unbounded ray of Problem \eqref{prob:s2}. When $\Obj_{\eqref{prob:s2}'} \le \Obj$, the term associated with $w$ in Problem \eqref{prob:independent:s2} can be disregarded, thus Problems \eqref{prob:independent:s2} and \eqref{prob:s2}$'$ have a finite optimum $\Obj_{\eqref{prob:independent:s2}}$ at $(\hat \psi, \hat u_y, 0)$. Otherwise, i.e., when Problem \eqref{prob:s2}$'$ is unbounded with an unbounded ray of $(\tilde \psi, \tilde u_y)$, Problem \eqref{prob:s2} is unbounded by $(\tilde \psi, \tilde u_y,0)$. \qed

\section{Proof of Corollary \ref{coro:independent2}.} 
Similar to the proof of Corollary \ref{coro:independent1}, it suffices to show that solving Problem \eqref{prob:s2}$''$ is sufficient to obtain the optimal solution or unbounded ray of Problem \eqref{prob:s2}. We define $(\psi',u_y') = (\frac{\psi}{w},\frac{u_y}{w})$, then Problem \eqref{prob:s2} becomes as follows:
	\begin{subequations}
	\begin{alignat}{2}
	\max \ & w \left[{\psi'}^T (b - A \hat x)+ {u_y'}^T (h_y - G_{xy} \hat x)-\Obj\right]  \\
    \mbox{ s.t.}  \ &B^T \psi' + G_y^T u'_y \preceq_{\K_y} d, \label{eq:independent2:s2:1}\\
     & \psi',u'_y, w \ge 0.
	\end{alignat}\label{prob:independent2:s2}
	\end{subequations}
The same analysis as in the proof of Corollary \ref{coro:independent1} holds; If Problem \eqref{prob:s2}$''$ has a finite optimum $\Obj_{\eqref{prob:s2}''}$ at $(\hat \psi, \hat u_y)$ and $\Obj_{\eqref{prob:s2}''} > \Obj$, $(\hat \psi, \hat u_y, 1)$ gives an unbounded ray of Problem \eqref{prob:s2}. If $\Obj_{\eqref{prob:s2}''} \le \Obj$, $(0, 0,0)$ is an optimal solution of Problem \eqref{prob:s2}. Otherwise, i.e., Problem \eqref{prob:s2}$''$ is unbounded by a feasible ray of $(\tilde \psi, \tilde u_y)$, Problem \eqref{prob:s2} is also unbounded by the feasible ray of $(\tilde \psi, \tilde u_y, 0)$. \qed

\section{Proof of Proposition \ref{prop:normalizing-feasibility-cut}.}
If $t(\lambda) < 0$, the optimal objective value of the inner
optimization problem of Problem \eqref{prob:LagRelax} 
approaches zero as $w$ converges to 0. If $t(\lambda) > 0$ then
the inner
optimization problem of Problem \eqref{prob:LagRelax} is unbounded. Therefore, Problem \eqref{prob:LagRelax} becomes equivalent to the following problem:
\begin{equation}
  \min_{\lambda \in \mathbb R} \left\{\lambda:  t(\lambda) \le 0 \right\}.
\label{prob:Lag2}
\end{equation} 
\noindent
Note that $t(\lambda)$ is non-increasing in $\lambda$. In addition, $(\hat y, \hat u_\psi, \hat v)$ and $(\frac{\tilde\psi}{\tilde w}, \frac{\tilde{u_y}}{\tilde{w}})$ are respectively feasible to Problems \eqref{prob:t1} and \eqref{prob:t2} when $\lambda = 0$, and thus $t(0) \ge
\frac{\U_\eqref{prob:s2}}{\tilde w} > 0$. Therefore, the optimal solution $\lambda^*$ of
Problem \eqref{prob:Lag2} is the solution of $t(\lambda) = 0.$ \qed

\section{Results on \texttt{MIPLIB} Instances}\label{appendix:result:miplib}
\begin{table}[H]
	\centering 
	\fontsize{8}{10}\selectfont
	  \caption{Computational Performance Comparison on \texttt{MIPLIB} Instances.\label{table:result:miplib}}
	  {\begin{tabular}[h]{crrrrrr}
		  \hline\noalign{\smallskip}
		   & \multicolumn{2}{c}{\cite{fischetti2017new}} & \multicolumn{2}{c}{\texttt{N}} & \multicolumn{2}{c}{\texttt{N+H}}\\
		  \noalign{\smallskip}\hline\noalign{\smallskip}
		  Instance &  Time (s) 	& Gap (\%)  &  Time (s) 	& Gap (\%) &  Time (s) 	& Gap (\%)\\
		  \noalign{\smallskip}\hline\noalign{\smallskip}   
		   air03-1    & $\dagger$ &   9.0 & $\dagger$ &  11.5 & $\dagger$ & 11.5 \\
	  \bf  air03-5    & $\dagger$ &	 29.5 &	$\dagger$ &	 35.9 & $\dagger$ & 26.1\\
	  \bf  air03-9    & $\dagger$ &  56.9 &	$\dagger$ &	 54.1 & $\dagger$ & 54.9\\
		   air04-1    & $\dagger$ &	  0.8 &	$\dagger$ &	  0.7 & $\dagger$ & 1.3\\
		   air04-5    & $\dagger$ &	  5.1 &	$\dagger$ &	  7.2 & $\dagger$ & 9.4\\
		   air04-9    & $\dagger$ &	 19.6 &	$\dagger$ &	 27.6 & $\dagger$ & 25.1\\
	  \bf  air05-1    & $\dagger$ &	  0.8 &	$\dagger$ &	  0.8 & $\dagger$ & 0.8\\
		   air05-5    & $\dagger$ &	 17.5 &	$\dagger$ &	 18.8 & $\dagger$ & 18.7\\
		   air05-9    & $\dagger$ &	 35.9 & $\dagger$ &  42.7 & $\dagger$ & 40.7\\
		   cap6000-1  &    478.27 &	  0.0 & $\dagger$ &  35.9 & $\dagger$ & 24.6\\
	  \bf  cap6000-5  & $\dagger$ & 106.9 & $\dagger$ &  91.0 & $\dagger$ & 50.0\\
	  \bf  cap6000-9  & $\dagger$ & 562.2 & $\dagger$ & 552.1 & $\dagger$ & 365.4\\
	  \bf  enigma-1   &      0.48 &   0.0 &      0.29 &   0.0 &      0.38 & 0.0\\
		   enigma-5   &     55.99 &   0.0 &    678.85 &   0.0 &    718.11 & 0.0\\
	  \bf  enigma-9   &      0.09 &   0.0 &      0.04 &   0.0 &      0.07 & 0.0\\
	  \bf  fast0507-1 &      7.92 &   0.0 &      0.43 &   0.0 &      0.56 & 0.0\\
	  \bf  fast0507-5 &      5.47 &   0.0 &      0.39 &   0.0 &      0.56 & 0.0\\
	  \bf  fast0507-9 &      3.78 &   0.0 &      0.47 &   0.0 &      0.67 & 0.0\\
	  \bf  harp2-1    & $\dagger$ &   5.0 & $\dagger$ &   7.3 & $\dagger$ & 1.7\\
	  \bf  harp2-5    & $\dagger$ &  25.1 & $\dagger$ &  68.4 & $\dagger$ & 18.2\\
		   harp2-9    & $\dagger$ & 130.5 & $\dagger\dagger$ &    -  & $\dagger\dagger$  & - \\
		   l152lav-1  &      8.35 &   0.0 &      3.36 &   0.0 &     30.78 & 0.0\\
		   l152lav-5  & $\dagger$ &   1.5 & $\dagger$ &  4.2 & $\dagger$ & 4.1\\
		   l152lav-9  & $\dagger$ &   5.5 & $\dagger$ & 6.5      & $\dagger$ & 6.5 \\
		   lseu-1     & 0.44      &   0.0 &      0.26 &   0.0 &      0.59 & 0.0\\
		   lseu-5     & $\dagger$ &  48.2 & $\dagger$ &  54.8 & $\dagger$ & 54.7 \\
	  \bf  lseu-9     & 0.97      &   0.0 &      0.34 &   0.0 &      0.52 & 0.0 \\
		   mitre-1    & $\dagger$ &   5.4 & $\dagger$ &   5.7 & $\dagger$ &  5.5 \\
	  \bf  mitre-5    & $\dagger$ &  21.5 & $\dagger$ &  21.7 & $\dagger$ & 21.2 \\
	  \bf  mitre-9    & $\dagger$ &  31.6 & $\dagger$ &  31.8 & $\dagger$ & 31.6 \\
		   mod010-1   & $\dagger$ &  0.04 & $\dagger$ &   0.15 & $\dagger$ & 0.07\\
		   mod010-5   & $\dagger$ &   2.0 & $\dagger$ &   3.9 & $\dagger$ & 3.9 \\
		   mod010-9   & $\dagger$ &  13.5 & $\dagger$ &  14.1 & $\dagger$ & 16.5 \\
	  \bf  nw04-1     & 906.94    &   0.0 &    348.19 &   0.0 &  488.62 & 0.0 \\
	  \bf  nw04-5     & $\dagger$ &  45.1 & $\dagger$ &  35.4 & $\dagger$ & 35.3 \\
	  \bf  nw04-9     & $\dagger$ &  65.7 & $\dagger$ &  59.8 & $\dagger$ & 59.5 \\
	  \bf  p0033-1    & 0.10      &   0.0 &      0.04 &   0.0 &      0.07 & 0.0\\
	  \bf  p0033-5    & 0.15      &   0.0 &      0.03 &   0.0 &      0.08 & 0.0 \\
		   p0033-9    & 0.04      &   0.0 &      0.02 &   0.0 &      0.05 &  0.0 \\
	  \bf  p0201-1    & $\dagger$ &  $\ddagger$ & $\dagger$ &  37.2 & $\dagger$ & 36.1\\
	  \bf  p0201-5    & $\dagger$ &  $\ddagger$ & $\dagger$ &  42.1 & $\dagger$ & 41.5 \\
		   p0201-9    & 1.07      &   0.0 &      0.22 &   0.0 &       9.95 & 0.0 \\
		   p0282-1    & $\dagger$ &   0.9 & $\dagger$ &   1.1 & $\dagger$ & 1.1\\
		   p0282-5    & $\dagger$ &   5.6 & $\dagger$ &   5.9 & $\dagger$ & 5.9 \\
		   p0282-9    & $\dagger$ &  33.5 & $\dagger$ &  40.3 & $\dagger$ & 39.8 \\
	  \bf  p0548-1    & $\dagger$ &  25.1 & $\dagger$ &  23.7 & $\dagger$ & 22.1\\
	  \bf  p0548-5    & $\dagger$ &  56.3 & $\dagger$ &  56.0 & $\dagger$ & 55.5 \\
		   p0548-9    & $\dagger$ &  36.3 & $\dagger$ &  39.9 & $\dagger$ & 39.6 \\
	  \bf  p2756-1    & $\dagger$ &  76.2 & $\dagger$ &  77.5 & $\dagger$ & 74.1\\
	  \bf  p2756-5    & $\dagger$ &  85.7 & $\dagger$ &  85.9 & $\dagger$ & 84.2 \\
	  \bf  p2756-9    & $\dagger$ &  88.4 & $\dagger$ &  88.5 & $\dagger$ & 87.8 \\
		   seymour-1  & $\dagger$ &   1.1 & $\dagger$ &   0.9 & $\dagger$ & 1.1\\
	  \bf  seymour-5  &      3.83 &   0.0 &      0.88 &   0.0 &      1.93 & 0.0 \\
	  \bf  seymour-9  &      0.31 &   0.0 &      0.06 &   0.0 &      0.05 & 0.0 \\
		   stein27-1  &      0.75 &   0.0 &      0.19 &   0.0 &      1.04 & 0.0\\
	  \bf  stein27-5  &      0.02 &   0.0 &      0.01 &   0.0 &      0.02 & 0.0 \\
	  \bf  stein27-9  &      0.01 &   0.0 &     0.006 &   0.0 &      0.01 & 0.0 \\
		   stein45-1  &      8.38 &   0.0 &      1.94 &   0.0 &      11.98 & 0.0\\
	  \bf  stein45-5  &      0.14 &   0.0 &      0.03 &   0.0 &      0.05 & 0.0 \\
		   stein45-9  &      0.01 &   0.0 &      0.005 &   0.0 &     0.02 & 0.0 \\
	\noalign{\smallskip}\hline
	  \end{tabular}}\\
	$\dagger$: The method times out; solution time $>$ 3,600 seconds);
	$\ddagger$: Numerical error occurs. The method is terminated with an infeasible solution;
	$\dagger\dagger$: Terminated with memory issues.
	\end{table}

	\section{Unit Commitment With Gas Awareness (UCGNA)}
\label{sec:UCGNA}


The UCGNA is a tri-level optimization problem where the first and second levels determine how to commit and dispatch electric power generating units; The third level decides how to operate the gas network given the natural gas demands of committed gas-fueled generators that are determined in the first and second levels. The economic feedback from the gas network, i.e., the natural gas zonal prices, is given by the dual solution $\psi$ of the third-level optimization and the first-level optimization is subject to constraints over both $\psi$ and commitment decisions $x$ in order to ensure the robustness of the unit commitment decisions against the economic feedback from the gas system. 

\cite{byeon2019unit} showed that the tri-level problem can be reformulated as a special case of BSOCP discussed in Section \ref{sec:sequence-independent}, which is in the form of Problem \eqref{prob:bl2} with $c_y = d$. The bilevel problem has a leader problem that decides the commitment decision (a subvector of $x$) and the follower problem is a joint network flow problem for dispatching electricity and natural gas with the given commitment decision $x$. Based on the follower’s dual solution $\psi$, which approximates the gas price, an additional constraint on both $x$ and $\psi$ is enforced in the leader problem to find a robust commitment decision against volatile natural gas prices in the gas system. 

\subsection{Mathematical Model}
This section specifies how the leader and the follower problem is formulated. In what follows, the electricity transmission grid is represented by an
undirected graph $\G^e = (\N, \E)$ and the natural gas transmission
system is by a directed graph $\G^g = (\V, \A)$. 
The letter $\T$ denotes
the set of time periods $\{0,1, \cdots, T\}$, and Tables \ref{table:param:e} and \ref{table:param:g}
summarize the parameters of the electricity and gas systems.
$[a,b]_{\mathbb Z}$ denotes the set of
integers in interval $[a,b]$, and $[n]$ denotes the set
$\{1,\cdots,n\}$ for some integer $n \ge 1$.

\begin{table}[!t]\fontsize{9}{9}\selectfont
	\centering
	\caption{Parameters of the Electricity System.} \label{table:param:e}
	{\begin{tabular}{p{0.1\textwidth} p{0.85\textwidth}}
			\toprule
			$\G^e = (\N, \E)$  & Undirected graph where $\N$ is a set of buses indexed by $i = 1, \cdots, N$ and $\E$ is a set of lines indexed with $l = 1, \cdots, E$ \\%
			$\Uc$ & Set of generators, indexed by $u = 1, \cdots, U$\\%
			$\quad\Uc^g \subseteq \Uc$ & Set of GFPPs \\
			$\quad\Uc(i) \subseteq \Uc$ & Set of generators located at $i \in \N$ \\%
			$\B_u$ & Set of supply bids submitted by $u \in \Uc$, indexed by $b = 1, \cdots, B_u$ \\
			$\quad c^e_{u,b}$ & Bid price of $b \in \B_u$\\ 
			$\quad\overline {s}_{u,b}$& Amount of 
			real power generation of $b \in \B_u$ \\%
			$\quad \mu_{u,b}$ & Maximum allowable gas price for bid $b$ to be profitable\\ 
			$\underline p_u, \overline p_u$ &  Minimum/maximum real power generation of $u \in \Uc$ \\%
			$\underline R_u, \overline R_u$ & Ramp-down/-up rate of $u \in \Uc$ \\
			$c_u$ & No-load cost of $u \in \Uc$\\%
			$\{H_{u,i}\}_{i = 0,1,2}$ & Coefficients of the heat rate curve of $u \in \Uc^g$\\%
			$\alpha_u$ & Maximum allowable percentage of the expense on natural gas over its marginal bid price for $u \in \Uc^g$\\
			$\Psi_{u}$ & Set of counts of time periods with distinct start-up costs of $u$ indexed by $h$ \\
			$\quad {C}_{u,h}$ & Start-up cost of $u\in \Uc$ when $u$ is turned on after it has been offline for some time $\in [\Psi_{u,h}, \Psi_{u,h+1}]$\\%
			$\overline o_{u,0}, \overline p_{u,0}$ & Initial on-off status/real power generation of $u \in \Uc$ \\%
			$ \underline\tau_{u}, \overline \tau_{u}$ &  Minimum-down/-up time of  $u \in \Uc$ \\  
			$\underline \tau_{u,0}, \overline \tau_{u,0}$ &  The time that generator $u \in \Uc$ has to be inactive/active from $t = 0$ \\%
			$\alpha_u$ & \\%
			$b_l$ & Line susceptance of $l \in \E$\\%
			$\overline f_{l}$ & Real power limit of $l \in \E$\\%
			$(d^e_{i,t})_{i \in \N}$ & Electricity load profile during $t \in \T$ \\%
			$\Delta_l$ & Maximum voltage angle difference between two end-points of $l \in \E$\\
			$\underline \theta_i, \overline \theta_i$ & Minimum/maximum voltage angle at $i \in \N$\\
			\bottomrule
	\end{tabular}}{}
\end{table}

\begin{table}[t!]\fontsize{9}{9}\selectfont
	\centering
	\caption{Parameters of the gas system} \label{table:param:g}
	{\begin{tabular}{p{0.1\textwidth} p{0.85\textwidth}}
			\toprule
			$\G^g = (\V, \A)$  & Directed graph representing a natural gas transmission network, where $\V$ is a set of junctions, indexed with $j = 1, \cdots, V$, and $\A \subseteq \V \times \V$ is a set of connections, indexed with $a = 1, \cdots, A$ \\%
			$\quad\A_c \subseteq \A$ & Set of compressors\\
			$\quad\A_v \subseteq \A$ & Set of control valves\\
			$\underline s^g_{j}, \overline s^g_{j}$ & Lower/Upper limit on natural gas supply at $j \in \V$\\%
			$\mathcal S_j$ & Set of non-overlapping intervals covering $[\underline s^g_j, \overline s^g_j]$, each with a distinct slope $c_{j,s}$ satisfying $c_{j,s} \le c_{j,s+1}$ for all consecutive intervals $s,s+1 \in \mathcal S_j$\\%
			$\kappa_j$ & Cost of demand shedding at $j \in \V$ \\%
			$ (d^g_{j,t})_{j \in \V}$ & Gas demand profile during $t \in \T$ \\%
			$W_a$ & Pipeline resistance (Weymouth) factor of $a \in \A$\\
			$\underline \pi_j, \overline \pi_j$ & Minimum/maximum squared pressure at $j \in \V$\\
			$\underline \alpha^c_a, \overline \alpha^c_a$ & Lower/upper compression ratio of $a \in \A_c$\\
			$\underline \alpha^v_a, \overline \alpha^v_a$ & Lower/upper control ratio of $a \in \A_v$\\
			$\mathcal K $ & Set of pricing zones, indexed with $k = 1, \cdots, K$\\
			$\quad \V(k) $ & Set of junctions that belong to $k \in \K$\\
			\bottomrule
	\end{tabular}}{}
\end{table}

\subsubsection{The Leader Problem}

The variables of the leader problem is summarized in Table \ref{table:var:l}. With
these notations, the leader model is specified in Problem \eqref{prob:l}. 
\begin{table}[!t]\fontsize{9}{9}\selectfont
	\centering
	\caption{Variables of the Leader Problem.}
        \label{table:var:l}
	{\begin{tabular}{p{0.05\textwidth} p{0.85\textwidth}}
			\toprule
			\multicolumn{2}{l}{  \textbf{Binary variables}} \\
			$o_{u,t}$ &  1 if $u \in \Uc$ is on during $t \in \T$, 0 otherwise \\
			$v^+_{u,t}$ &  1 if $u \in \Uc$ becomes online during $t \in \T$, 0 otherwise \\
			$v^-_{u,t}$ &  1 if $u \in \Uc$ becomes offline during $t \in \T$, 0 otherwise \\
			$w_{u,b,t}$ &  1 if $b \in \B_u$ is selected during $t \in \T$, 0 otherwise\\
			\multicolumn{2}{l}{  \textbf{Continuous variables}} \\
			$r_{u,t}$ & Start-up cost of $u \in \Uc$ during $t \in \T$ \\
			$\varphi_{u,t}$& Maximum allowable natural gas price for $u \in \Uc^g$ to generate power at its scheduled level during $t \in \T$\\
			\bottomrule
	\end{tabular}}{}\end{table}

\begin{figure}[!t]
\begin{subequations}\fontsize{9}{9}\selectfont
    \begin{equation}
\min \  \sum_{t\in [T]}\left( \beta \sum_{u \in \Uc}  (
		 c_u o_{u,t} +  r_{u,t}) + \left(\beta \sum_{u \in \Uc} \sum_{b \in \B_u}  c^e_{u,b} {s}^e_{u,b,t}  + (1-\beta) \sum_{j \in \V}  (\sum_{s \in \mathcal S_j}  c^g_{j,s} s^g_{s,t} + \kappa_j q_{j,t})\right) \right)\label{e:obj}    
    \end{equation}
	\begin{align}
		\mbox{s.t.} \  
		&r_{u,t} \ge C_{u,h} ( o_{u,t} - \sum_{n \in [h]} o_{u,t-n} ),  &  \forall h \in \Psi_s, u \in \Uc, t \in [T], \label{e:UC:su}\\
		&r_{u,t} \ge 0, & \forall u \in \Uc, t \in [T], \label{e:UC:su:nonnegative}\\
		&o_{u,t} = \overline o_{u,0}, & \forall u \in \Uc, \  t \in [0,\overline\tau_{u,0} + \underline \tau_{u,0}]_{\mathbb Z}, \label{e:UC:initial_gen_status}\\
		& \sum_{t' \in [t - \overline\tau_{u} +1, t]_{\mathbb Z}} v^+_{u,t'} \le o_{u,t},  & \forall  u \in \Uc, \ t \in [\max\{\overline\tau_{u}, \overline\tau_{u,0} + 1\}, T]_\mathbb{Z}, \label{e:UC:min_up}\\
		& \sum_{t' \in [t - \underline \tau_{u} +1,t]_{\mathbb Z} } v^+_{u,t'} \le 1- o_{u,t - \underline \tau_{u}},  & \forall u \in \Uc, \ t \in [\max\{\underline \tau_{u}, \underline \tau_{u,0} + 1\}, T]_\mathbb{Z}, \ \label{e:UC:min_down}\\
		& v^+_{u,t} - v^-_{u,t} = o_{u,t} - o_{u,t-1}, &\forall u \in \Uc, \  t \in [T], \label{e:UC:logic:on_off}\\
		&  w_{u,b,t} \le o_{u,t}, & \forall b \in \B_u, \ u \in \Uc^g, \ t \in [T], \label{e:UC:logic:bid1}\\
        &\varphi_{u,t} = \sum_{b \in [B_u-1]} \mu_{u,b} (w_{b,t} - w_{u,b+1,t}) + \mu_{u,B_u} w_{u,B_u,t}, & \forall u \in \Uc^g, \ t \in [T], \label{e:bid_price}\\
		& v^+_{u,t}, v^-_{u,t}, o_{u,t} \in \{0,1\}, & \forall u \in \Uc, \ t \in [T], \label{e:UC:logic:nonnega}\\
		& w_{u,b,t} \in \{0,1\}, & \forall b \in \B_u, u \in \Uc^g, \ t \in [T], \label{e:UC:logic:nonnega:w}\\
		&  0 \le s^e_{u,b,t} \le \overline {s}_{u,b}w_{u,b,t}, & \forall b \in \B_u, \ u \in \Uc^g, \ t \in [T], \label{e:UC:bound:bids:gfpp}\\
		&  \overline s_{u,b} w_{u,b+1,t} \le s^e_{u,b,t}, & \forall b \in [1, B_u-1]_{\mathbb Z}, \ u \in \Uc^g, \ t \in [T], \label{e:UC:logic:bid2}\\
		& \varphi_{u,t} \ge \psi_{k,t} o_{u,t},& \forall k \in \K,  i \in  \V(k),  u \in \U(i) \cap \Uc^g, \ t \in [T], \label{c:bid_val}\\
		& (y,\psi) \in \Q \left(\mbox{Problem } \eqref{prob:f}\right).\label{l:f}
	\end{align}
	\label{prob:l}
\end{subequations}
\end{figure}

The objective function includes the objective of the unit-commitment problem (i.e., the no-load
costs, the start-up costs, and the costs of the selected supply bids
of each electrical power generating units $\sum_{t \in [T]} \sum_{u \in \Uc}\left( c_u o_{u,t} +  r_{u,t}+ \sum_{b \in \B_u}  c^e_{u,b} {s}^e_{u,b,t}\right)$) and the cost of dispatching natural gas ($\sum_{t \in [T]}\sum_{j \in \V}  (\sum_{s \in \mathcal S_j}  c^g_{j,s} s^g_{s,t} + \kappa_j q_{j,t})$), which are respectively scaled by $\beta \in (0,1)$ and $1-\beta$. Equation \eqref{e:UC:su}
computes the start-up cost $r_{u,t}$ of a generator $u$ for time
period $t$ based on how long $u$ has been offline \citep{morales2013tight}.  The
expression $o_{u,t} - \sum^h_{n=1} o_{u,t-n}$ is one when generator
$u$ becomes online after it has been turned off for $h$ time
periods. Equation \eqref{e:UC:su:nonnegative} states the nonnegativity
requirement on $r_{u,t}$. Equation \eqref{e:UC:initial_gen_status}
specifies the initial on-off status of each generator. The minimum-up
and -down constraints are specified in Equations \eqref{e:UC:min_up}
and \eqref{e:UC:min_down} respectively. The relationship between the
variables for the on-off, start-up, and shut-down statuses of each
generator is stated in Equation \eqref{e:UC:logic:on_off}. 
Equation \eqref{e:UC:logic:bid1}
states that the bid of a generator can be selected only when it is
committed. Equation
\eqref{e:UC:bound:bids:gfpp} is bound constraints for the bids
submitted by the GFPPs, which ensures that the
indicator variable $w_{b,t}$ is one whenever bid $b$ is used for time
period $t$ (i.e., $s^e_{b,t} > 0$). In Equation \eqref{e:UC:logic:bid2}, the
$(b+1)^{\mbox{th}}$ bid is selected only if the bid $b$ is fully
used. Accordingly, Equation \eqref{e:bid_price} states that
$\varphi_{u,t}$ is the maximum allowable gas price for $u \in \Uc^g$
to be profitable when generating its scheduled amount. The binary
requirements for logical variables $v^+_{u,t}, v^-_{u,t},o_{u,t}$ and $w_{b,t}$ are
specified in Equations \eqref{e:UC:logic:nonnega} and \eqref{e:UC:logic:nonnega:w}.
The economic coupling between the electricity and gas networks is
enforced by {\em bid-validity constraints} (Equation \eqref{c:bid_val}) that state that the power generation of a committed gas-fired power plant $u \in \Uc^g$, which receives natural gas at junction $k \in \V$, should be profitable with regard to the realized natural gas price $\psi_{k,t}$. The nonlinear term in the right-hand side of Equation \eqref{c:bid_val} is linearized by employing an exact McCormick relaxation: For each $k \in \K,  i \in  \V(k),  u \in \U(i) \cap \Uc^g, \ t \in [T],$
\begin{subequations}
\begin{align}
    \varphi_{u,t} \ge \upsilon_{u,k,t},\\
    \upsilon_{u,k,t} \ge \psi_{k,t} -\overline{\psi}_{k,t} (1-o_{u,t}),\\
    \upsilon_{u,k,t} \le \psi_{k,t} -\underline{\psi}_{k,t} (1-o_{u,t}),\\
    \upsilon_{u,k,t} \le \overline{\psi}o_{u,t},\\
    \upsilon_{u,k,t} \ge \underline{\psi}o_{u,t}.
\end{align}
\label{eqs:mcCormick}
\end{subequations}
Although the natural gas system is operated in a
decentralized manner, the zonal price of natural gas $\boldsymbol\psi$
can be approximated with the dual solution of the follower problem. In Equation \eqref{l:f}, $\mathcal Q$ denotes the projection of optimal pairs of primal and dual solutions of the follower problem (Problem \eqref{prob:f}) onto the space of $\boldsymbol{s}^e$ and $\boldsymbol \psi$. The bid validity constraints use the maximum natural gas price (e.g.,
\$200 per mmBtu) as $\overline \psi$ and 0 as $\underline \psi$.

Note that Equations \eqref{e:UC:su}-\eqref{e:UC:logic:nonnega:w} are the specification of $\X$ in Equation \eqref{eq:bl2:upper:x} and Equations \eqref{e:UC:bound:bids:gfpp} and \eqref{e:UC:logic:bid2} are that of Equation \eqref{eq:bl2:upper:y}. Equation \eqref{eq:bl2:upper:psi} is specified by Equations \eqref{eqs:mcCormick}.

\subsubsection{The Follower Problem}
\begin{table}[!t]\fontsize{9}{9}\selectfont
	\centering
	\caption{Variables of the Follower Problem.}
        \label{table:var:f}
	{\begin{tabular}{p{0.05\textwidth} p{0.7\textwidth}}
			\toprule
			\multicolumn{2}{l}{  \textbf{Variables on the electricity system}} \\
			$s^e_{b,t}$ & Real power generation from $b \in \B_u$ of $u \in \Uc$ during $t \in \T$\\
			$p_{u,t} $ & Real power generation of $u \in \Uc$ during $t \in \T$ \\
			$f_{l,t}$ & Real power flow on $l \in \E$ during $t \in \T$ \\
			$\theta_{i,t}$ & Voltage angle on $i \in \N$ during $t \in \T$\\
			\multicolumn{2}{l}{ \textbf{Variables on the gas system}} \\
			$s^g_{j,t}$ & Amount of gas supplied at $j \in \V$ during $t \in \T$ \\
			$s^g_{j,s,t}$ & Amount of gas supply from $s \in \mathcal S_j$ during $t \in \T$\\
			$\pi_{j,t}$ & Pressure squared at $j \in \V$ during $t \in \T$\\
			$\phi_{a,t}$ & Gas flow on $a \in \A$ during $t \in \T$\\
			$l_{j,t}$ & Satisfied gas demand at $j \in \V$ during $t \in \T$\\
			$q_{j,t}$ & Shedded gas demand at $j \in \V$ during $t \in \T$\\
			$\gamma_{j,t}$ & Total amount of gas consumed by the GFPP located at $j \in \N \cap \V$ during $t \in \T$\\
			\bottomrule
	\end{tabular}}{}\end{table}
	
\begin{figure}[!t]
\begin{subequations}\fontsize{9}{9}\selectfont
  \begin{align}
		\min \quad &  \sum_{t \in [T]}\left(\beta \sum_{u \in \U} \sum_{b \in \B_u}  c^e_{u,b}s^e_{u,b,t}  + (1-\beta)  \sum_{j \in \V}  (\sum_{s \in \mathcal S_j}  c_{j,s} s^g_{j,s,t} + \kappa_j q_{j,t})\right) \label{e:ED:obj}\\
		  \mbox{s.t.} \quad & \sum_{u \in \U(i)} {p}_{u,t} - d^e_{i,t}= \sum_{l \in \E: l_t = i} f_{l,t} - \sum_{l \in \E: l_h = i} f_{l,t}, & \forall i \in \N,\ t \in [T], \label{e:ED:bal}\\
		  & {p}_{u,t} = \sum_{b \in \B_u} {s}^e_{b,t}  & \forall u \in \Uc, \ t \in [T], \label{e:ED:bid}\\
		  &  0 \le s^e_{b,t} \le \overline {s}_{b},  & \forall b \in \B_u, \ u \in \Uc, \ t \in [T], \label{e:ED:bound:bids} \\
		  &  \underline{p}_{u}o_{u,t} \le p_{u,t} \le \overline{p}_{u} o_{u,t},  & \forall u \in \Uc, \ t \in [T], \label{e:ED:bound:gen} \\
		  &p_{u,0} = \overline p_{u,0}, & \forall u \in \Uc, \label{e:ED:initial_gen}\\
		  & p_{u,t} - p_{u,t-1} \le \overline R_u o_{u,t-1} + \overline p_{u} v^+_{u,t}, & \forall   u \in \Uc, \ t \in [T], \label{e:ED:rup}\\
		& p_{u,t-1} - p_{u,t} \le \underline R_u o_{u,t-1} + \underline p_{u} v^-_{u,t}, & \forall   u \in \Uc, \ t \in [T], \label{e:ED:rdown}\\
		& f_{l,t} = -b_l (\theta_{l_h,t} - \theta_{l_t,t}),  & \forall l \in \E, \ t \in [T], \label{e:ED:DCOPF:PF}\\
		& -\overline f_l \le f_{l,t} \le \overline f_l, & \forall l \in \E,  \ t \in [T], \label{e:ED:DCOPF:thermalLimit}\\
		& \underline \theta_{i} \le \theta_{i,t}\le \overline \theta_{i}, & \forall i \in \N, \ t \in [T], \label{e:ED:DCOPF:angleBound}\\
		& -\Delta_{l}  \le \theta_{l_h,t} - \theta_{l_t,t} \le \Delta_{l}  & \forall l \in \E, \ t \in [T], \label{e:ED:DCOPF:angleDiff}\\
		& s^g_{j,t} - l_{j,t} - \gamma_{j,t} = \sum_{a \in \A: a_t = j} \phi_{a,t} - \sum_{a \in \A: a_h = j} \phi_{a,t}, & \forall j \in \V, t \in [T], \label{g:bal}\\
		& s^g_{j,t} = \sum_{s \in \mathcal S_j} s^g_{j,s,t}, & \forall j \in \V, \ t \in [T], \label{g:link}\\
		& l_{j,t} = d^g_{j,t} - q_{j,t}, & \forall j \in \V, t \in [T],\label{g:load_shed}\\
		& 0 \le q_{j,t} \le d^g_{j,t}, & \forall j \in \V, t \in [T],\label{g:load_shed_bound}\\
		& \phi_{a,t} \ge 0, & \forall a \in \A, t \in [T],\label{g:flux_bound}\\
		&  \underline s^g_{j} \le s^g_{j,t} \le \overline s^g_{j},& \forall j \in \V,\ t \in [T],\label{g:bound:supply}\\
		& \underline \alpha^c_a \pi_{a_h,t} \le \pi_{a_t,t}  \le  \overline \alpha^c_a \pi_{a_h,t}, & \forall a  \in \A_c,\ t \in [T], \label{g:compressors}\\
		& \underline \alpha^v_a \pi_{a_h, t}  \le \pi_{a_t,t}  \le  \overline \alpha^v_a \pi_{a_h,t}, & \forall a \in \A_v,\ t \in [T], \label{g:valves}\\
		&  \pi_{a_h, t}  - \pi_{a_t,t}  \ge W_a \phi_{a,t}^2,& \forall a \in \A \setminus (\A_v \cup \A_c),\ t \in [T], \label{g:Weymouth}\\
		&\underline \pi_j  \le \pi_{j, t}  \le \overline \pi_{j},& \forall j \in \V,\ t \in [T] \label{g:bound:pi} \\
        & \gamma_{j,t} \ge \sum_{u \in \U(i) \cap \Uc^g} H_{u,2} p_{u,t}^2 + H_{u,1} p_{u,t} + H_{u,0}, & \forall j \in \N \cap \V, \ t \in [T].	\label{c:physical}
		\end{align}
		\label{prob:f}
\end{subequations}
\end{figure}

Based on the commitment decisions decided in the leader problem, the follower problem (i.e.,
Equations \eqref{e:ED:obj} - \eqref{c:physical}) decides the
hourly operating schedule of each committed generators and the gas transmission network in order to
minimize the system costs for electricity ($\sum_{t \in [T]} \sum_{u \in \U} \sum_{b \in \B_u}  c^e_{u,b} {s}^e_{u,b,t}$) and gas ($\sum_{t \in [T]} \sum_{j \in \V} (\sum_{s \in \mathcal S_j}  c_{j,s} s^g_{j,s,t} + \kappa_j q_{j,t})$) which are respectively scaled by $\beta \in (0,1)$ and $1-\beta$. Equation \eqref{e:ED:bal} states
the flow conservation constraints for real power at each bus, using
$l_h$ and $l_t$ to represent the head and tail of $l \in \E$. Equation
\eqref{e:ED:bid} states that the total real power generation of a
generator $u$ is equal to the production of its selected bids.
Equation \eqref{e:ED:bound:bids} constrains the power generation
$s_{b,t}^e$ from bid $b \in \B_u$ to be no more than the submitted amount
$\bar s_b$. Equation \eqref{e:ED:bound:gen} enforces the bound on the
real power generation of each generator. Equation
\eqref{e:ED:initial_gen} specifies the initial generation amount of
each generator, and Equations \eqref{e:ED:rup} and \eqref{e:ED:rdown}
state the ramp-up and -down constraints of each generator. Equation
\eqref{e:ED:DCOPF:PF} captures the DC approximation of the power flow
equations and Equation \eqref{e:ED:DCOPF:thermalLimit} specifies the
thermal limit on each line. Equations \eqref{e:ED:DCOPF:angleBound}
and \eqref{e:ED:DCOPF:angleDiff} state the voltage angle bounds on
each bus and the bounds on the angle difference of two adjacent buses
respectively.

A steady-state natural gas model is specified in Equations \eqref{g:bal}-\eqref{g:bound:pi}, which is similar to those in
\cite{bent2018joint,sanchez2016convex,borraz2016convex} and uses the
Weymouth equation to capture the relationship between pressures and
flux. The flux conservation constraint is given in Equation
\eqref{g:bal}, where $a_h$ and $a_t$ represent the head and tail of $a
\in \A$. Equation \eqref{g:link} calculates the total gas production at junction $j \in \V$, and Equation \eqref{g:load_shed} determines the demand served at
each junction: It captures the amount of gas load shedding which must
be nonnegative and cannot exceed the demand at the corresponding
junction (Equation \eqref{g:load_shed_bound}). The model assumes that
gas flow directions are predetermined and Equation
\eqref{g:flux_bound} enforces the sign of gas flow variables, i.e., it
constrains $\phi_{a,t}$ to be nonnegative. Equation
\eqref{g:bound:supply} specifies the upper and lower limits of natural
gas supplies. The change in pressure through compressors and control
valves are formulated in Equations \eqref{g:compressors} and
\eqref{g:valves} and the model use a single compressor machine
approximation as in prior work.  The steady-state physics of gas flows
is formulated with the Weymouth equation in Equation
\eqref{g:Weymouth}. Equation \eqref{g:bound:pi} states the bounds on
nodal pressures. Equation \eqref{g:Weymouth} is a second-order cone relaxation of the Weymouth equation ($\pi_{a_h, t} - \pi_{a_t,t} = W_a \phi_{a,t}^2$) from \cite{borraz2016convex}, and the result therein empirically showed the relaxation is very tight. 

Gas-fired power plants also play as a physical interface between the electrical power and gas networks. The real power generation $\boldsymbol{p}$ of a gas-fired power plant induces a demand
$\boldsymbol{\gamma}$ in the natural gas system. Equation
\eqref{c:physical} specifies the relationship between the real power
generation of a gas-fueled generator and the amount of natural gas needed for the
generation. In the equation, this relationship is approximated by a
quadratic heat-rate curve ($\gamma_{j,t} = \sum_{u \in \U(i) \cap \Uc^g} H_{u,2} p_{u,t}^2 + H_{u,1} p_{u,t} + H_{u,0},  \forall j \in \N \cap \V, \ t \in [T]$), whose coefficients are given as $H_u$. The
equation is convexified like the Weymouth equation in Equation \eqref{c:physical}.

\section{Results on UCGNA Instances} \label{appendix:result:UCGNA}
	\begin{table}[!h]
		\centering \fontsize{9}{10}\selectfont
			\caption{Computational Performance Comparison ($\eta_p = 1$).\label{table:computation1}}
			{\begin{tabular}[h]{cccccccc}
					\hline\noalign{\smallskip}
					\multicolumn{2}{c}{Instance} & \multicolumn{2}{c}{\texttt{D}} & \multicolumn{2}{c}{\texttt{G}} & \multicolumn{2}{c}{\texttt{B}}\\
					\noalign{\smallskip}\hline\noalign{\smallskip}
					 $\eta_p$& $\eta_g$& Time (s) 	& Gap (\%)  &  Time (s) 	& Gap (\%) &  Time (s) 	& Gap (\%)\\
					\noalign{\smallskip}\hline\noalign{\smallskip}  
					\multirow{14}{*}{1} 
					 & 1 &  25.42 &   0.0 &   15.28   &   0.0 &$\dagger$ & 6.8\\
					 & 1.1 &  25.91&	0.0 &	23.24&	0.0 & $\dagger$ &4.3\\
					 & 1.2 &  25.86	&0.0	&	14.78&	0.0 & $\dagger$ & 2.2\\
					 & 1.3 &  29.33&	0.0	&	31.17&	0.0 & $\dagger$ & 4.4\\
					 & 1.4 &  26.60&	0.0	&	6.76&	0.0 & $\dagger$ & 2.6\\
					 & 1.5 &  25.80&	0.0&		13.24&	0.0 & $\dagger$ &6.2\\
					 & 1.6 &  27.01&	0.0	&	33.56&	0.0 & $\dagger$ & 3.1\\
					 & 1.7 &  100.82&	0.0&		22.78&	0.0 & $\dagger$ & 4.5\\
					 & 1.8 &  $\dagger$ &	1.8 &  	$\dagger$ &$\ddagger$ & $\dagger$ & 6.7\\
					 & 1.9 &  $\dagger$&	1.3 &  $\dagger$ & $\ddagger$ & $\dagger$ & 10.6\\
					 & 2.0 &  67.13 &   0.0 &    $\dagger$ &   1.3 & $\dagger$ & 10.8\\
					 & 2.1 &  1091.88   &   0.0 &   $\dagger$ &   3.2 & $\dagger$ & 20.0\\
					 & 2.2 &  566.94    &   0.0    &   $\dagger$    &   3.6 & $\dagger$ &19.1\\
					 & 2.3 &  31.52 &   0.0 &   15.94   &   0.0 & $\dagger$ & 8.4\\
		\noalign{\smallskip}\hline
			\end{tabular}}
		\end{table}
	\begin{table}[!t]
	\centering \fontsize{9}{10}\selectfont
		\caption{Computational Performance Comparison ($\eta_p = 1.3$).\label{table:computation13}}
		{\begin{tabular}[h]{cccccccc}
				\hline\noalign{\smallskip}
				\multicolumn{2}{c}{Instance} & \multicolumn{2}{c}{\texttt{D}} & \multicolumn{2}{c}{\texttt{G}} & \multicolumn{2}{c}{\texttt{B}}\\
				\noalign{\smallskip}\hline\noalign{\smallskip}
				$\eta_p$& $\eta_g$& Time (s) 	& Gap (\%)  &  Time (s) 	& Gap (\%) &  Time (s) 	& Gap (\%)\\
				\noalign{\smallskip}\hline\noalign{\smallskip}  
				  \multirow{14}{*}{1.3} 
				&	1	&	31.01	&	0.0	&	4.37	&	0.0	&	$\dagger$ &1.9\\
	&	1.1	&	28.93	&	0.0	&	3.20	&	0.0	&	$\dagger$&2.8\\
	&	1.2	&	30.87	&	0.0	&	3.28	&	0.0	&	$\dagger$&2.9\\
	&	1.3	&	48.22	&	0.0	&	2.93	&	0.0	&	$\dagger$&3.3\\
	&	1.4	&	32.69	&	0.0	&	12.07	&	0.0	&	$\dagger$&3.8\\
	&	1.5	&	44.13	&	0.0	&	23.89	&	0.0	&	$\dagger$&2.2\\
	&	1.6	&	$\dagger$	&	0.3	&	$\dagger$	&	0.2	&	$\dagger$&4.1\\
	&	1.7	&	$\dagger$	&	3.5	&	$\dagger$	&	$\ddagger$	&	$\dagger$&11.0\\
	&	1.8	&	$\dagger$	&	3.2	&	$\dagger$	&	$\ddagger$	&	$\dagger$&10.9\\
	&	1.9	&	$\dagger$	&	3.3	&	$\dagger$	&	$\ddagger$	&	$\dagger$&17.4\\
	&	2	&	$\dagger$	&	4.2	&	$\dagger$	&	19.9	&	$\dagger$&14.9\\
	&	2.1	&	$\dagger$	&	4.3	&	$\dagger$	&	$\ddagger$	&	$\dagger$&9.7\\
	&	2.2	&	$\dagger$	&	4.0	&	$\dagger$	&	$\ddagger$	&	$\dagger$&14.8\\
	&	2.3	&	43.23	&	0.0	&	10.43	&	0.0	&	$\dagger$&5.7\\
	\noalign{\smallskip}\hline
		\end{tabular}}
	\end{table}
	\begin{table}[!t]
	\centering \fontsize{9}{10}\selectfont
		\caption{Computational Performance Comparison ($\eta_p = 1.6$).\label{table:computation16}}
		{\begin{tabular}[h]{cccccccc}
				\hline\noalign{\smallskip}
				\multicolumn{2}{c}{Instance} & \multicolumn{2}{c}{\texttt{D}} & \multicolumn{2}{c}{\texttt{G}} & \multicolumn{2}{c}{\texttt{B}}\\
				\noalign{\smallskip}\hline\noalign{\smallskip}
				$\eta_p$& $\eta_g$& Time (s) 	& Gap (\%)  &  Time (s) 	& Gap (\%) &  Time (s) 	& Gap (\%)\\
			\noalign{\smallskip}	\hline
			\noalign{\smallskip}  \multirow{14}{*}{1.6} 
			&	1	&	43.51	&	0.0	&	4.33	&	0.0	&	$\dagger$	&5.8\\
	&	1.1	&	27.88	&	0.0	&	5.46	&	0.0	&	$\dagger$	&2.8\\
	&	1.2	&	26.63	&	0.0	&	7.67	&	0.0	&	$\dagger$	&3.9\\
	&	1.3	&	22.19	&	0.0	&	6.25	&	0.0	&	$\dagger$	&2.7\\
	&	1.4	&	29.75	&	0.0	&	6.35	&	0.0	&	$\dagger$	&4.7\\
	&	1.5	&	330.88	&	0.0	&	21.08	&	0.0	&	$\dagger$	&7.0\\
	&	1.6	&	$\dagger$	&	2.1	&	$\dagger$	&	$\ddagger$	&	$\dagger$	&9.7\\
	&	1.7	&	$\dagger$	&	2.0	&	$\dagger$	&	$\ddagger$	&	$\dagger$	&8.1\\
	&	1.8	&	$\dagger$	&	6.2	&	$\dagger$	&	$\ddagger$	&	$\dagger$	&17.1\\
	&	1.9	&	$\dagger$	&	7.4	&	$\dagger$	&	$\ddagger$	&	$\dagger$	&11.5\\
	&	2	&	$\dagger$	&	3.7	&	$\dagger$	&	$\ddagger$	&	$\dagger$	&8.7\\
	&	2.1	&	$\dagger$	&	5.0	&	$\dagger$	&	$\ddagger$	&	$\dagger$	&9.1\\
	&	2.2	&	$\dagger$	&	5.0	&	$\dagger$	&	$\ddagger$	&	$\dagger$	&9.0\\
	&	2.3	&	12.44	&	0.0	&	3.76	&	0.0	&	$\dagger$	&3.9\\
	\noalign{\smallskip}\hline
		\end{tabular}}
	\end{table}
	\begin{table}[!t]
	\centering \fontsize{9}{10}\selectfont
		\caption{Benefits of the Acceleration Schemes ($\eta_p = 1$).\label{table:accel1}}
		{\begin{tabular}[h]{cccccccccc}
			\hline\noalign{\smallskip}
				 & \multicolumn{2}{c}{\texttt{D($n_1, i_1$)}} & \multicolumn{2}{c}{\texttt{D($n_0, i_1$)}} & 
				\multicolumn{2}{c}{\texttt{D($n_1, i_0$)}} &
				\multicolumn{2}{c}{\texttt{D($n_0, i_0$)}} \\
				\noalign{\smallskip}\hline\noalign{\smallskip}
				  $\eta_g$& Time (s)	& Gap (\%) &  Time (s)	& Gap (\%) &  Time (s)	& Gap (\%) & Time (s)	& Gap (\%) \\
				\noalign{\smallskip}\hline\noalign{\smallskip}  
		1	&	25.42	&	0.00	&	30.61	&	0.00	&	$\dagger$	&	49.38	&	$\dagger$	&	52.65	\\
		1.1	&	25.91	&	0.00	&	25.39	&	0.00	&	$\dagger$	&	50.70	&	$\dagger$	&	52.53	\\
		1.2	&	25.86	&	0.00	&	25.35	&	0.00	&	$\dagger$	&	51.13	&	$\dagger$	&	53.59	\\
		1.3	&	29.33	&	0.00	&	28.19	&	0.00	&	$\dagger$	&	50.82	&	$\dagger$	&	52.67	\\
		1.4	&	26.60	&	0.00	&	26.74	&	0.00	&	$\dagger$	&	53.15	&	$\dagger$	&	53.20	\\
		1.5	&	25.80	&	0.00	&	27.51	&	0.00	&	$\dagger$	&	51.99	&	$\dagger$	&	52.63	\\
		1.6	&	27.01	&	0.00	&	25.90	&	0.00	&	$\dagger$	&	38.88	&	$\dagger$	&	53.36	\\
		1.7	&	100.82	&	0.00	&	98.52	&	0.00	&	$\dagger$	&	19.33	&	$\dagger$	&	53.30	\\
		1.8	&	$\dagger$	&	1.77	&	$\dagger$	&	1.42	&	$\dagger$	&	3.09	&	$\dagger$	&	52.81	\\
		1.9	&	$\dagger$	&	1.32	&	$\dagger$	&	1.47	&	$\dagger$	&	1.52	&	$\dagger$	&	53.36	\\
		2	&	67.13	&	0.00	&	58.85	&	0.00	&	$\dagger$	&	9.17	&	$\dagger$	&	52.96	\\
		2.1	&	1091.88	&	0.00	&	$\dagger$	&	4.80	&	$\dagger$	&	4.52	&	$\dagger$	&	52.56	\\
		2.2	&	566.94	&	0.00	&	$\dagger$	&	4.45	&	$\dagger$	&	5.23	&	$\dagger$	&	53.46	\\
		2.3	&	31.52	&	0.00	&	23.85	&	0.00	&	$\dagger$	&	38.59	&	$\dagger$	&	52.97	\\
				 \noalign{\smallskip}			\hline
		\end{tabular}}
	\end{table}
	\begin{table}[!t]
	\centering\fontsize{9}{10}\selectfont
		\caption{Benefits of the Acceleration Schemes ($\eta_p = 1.3$).\label{table:accel13}}
		{\begin{tabular}[!t]{ccccccccc}
				\hline\noalign{\smallskip}
				 & \multicolumn{2}{c}{\texttt{D($n_1, i_1$)}} & \multicolumn{2}{c}{\texttt{D($n_0, i_1$)}} & 
				\multicolumn{2}{c}{\texttt{D($n_1, i_0$)}} &
				\multicolumn{2}{c}{\texttt{D($n_0, i_0$)}} \\
				\noalign{\smallskip}\hline\noalign{\smallskip}
				  $\eta_g$& Time (s)	& Gap (\%) &  Time (s)	& Gap (\%) &  Time (s)	& Gap (\%) & Time (s)	& Gap (\%) \\
				\noalign{\smallskip}\hline\noalign{\smallskip}  
		1	&	31.01	&	0.00	&	30.83	&	0.00	&	$\dagger$	&	63.96	&	$\dagger$	&	63.78	\\
		1.1	&	28.93	&	0.00	&	27.83	&	0.00	&	$\dagger$	&	54.30	&	$\dagger$	&	63.93	\\
		1.2	&	30.87	&	0.00	&	143.36	&	0.00	&	$\dagger$	&	60.95	&	$\dagger$	&	63.65	\\
		1.3	&	48.22	&	0.00	&	52.89	&	0.00	&	$\dagger$&	56.01	&	$\dagger$&	64.09	\\
		1.4	&	32.69	&	0.00	&	31.04	&	0.00	&	$\dagger$	&	51.67	&	$\dagger$	&	64.85	\\
		1.5	&	44.13	&	0.00	&	44.98	&	0.00	&	$\dagger$	&	53.98	&	$\dagger$	&	64.80	\\
		1.6	&	$\dagger$	&	0.31	&	$\dagger$	&	1.08	&	$\dagger$	&	1.94	&	$\dagger$	&	65.07	\\
		1.7	&	$\dagger$	&	3.53	&	$\dagger$	&	5.34	&	$\dagger$	&	3.42	&	$\dagger$	&	65.99	\\
		1.8	&	$\dagger$	&	3.15	&	$\dagger$	&	4.01	&	$\dagger$	&	3.73	&	$\dagger$	&	65.92	\\
		1.9	&	$\dagger$	&	3.26	&	$\dagger$	&	8.28	&	$\dagger$	&	7.97	&	$\dagger$	&	66.22	\\
		2	&	$\dagger$	&	4.24	&	$\dagger$	&	4.59	&	$\dagger$	&	4.51	&	$\dagger$	&	64.58	\\
		2.1	&$\dagger$	&	4.27	&	$\dagger$	&	4.12	&	$\dagger$	&	4.29	&	$\dagger$	&	63.36	\\
		2.2	&	$\dagger$	&	4.03	&	$\dagger$	&	4.07	&	$\dagger$	&	4.08	&	$\dagger$	&	64.46	\\
		2.3	&43.23	&	0.00	&	48.06	&	0.00	&	$\dagger$	&	14.51	&$\dagger$&	62.93	\\
			\noalign{\smallskip}	\hline
		\end{tabular}}
	\end{table}
	\begin{table}[!t]
	\centering \fontsize{9}{10}\selectfont
		\caption{Benefits of the Acceleration Schemes ($\eta_p = 1.6$).\label{table:accel16}}
		{\begin{tabular}[!t]{cccccccccc}
				\hline\noalign{\smallskip}
				 & \multicolumn{2}{c}{\texttt{D($n_1, i_1$)}} & \multicolumn{2}{c}{\texttt{D($n_0, i_1$)}} & 
				\multicolumn{2}{c}{\texttt{D($n_1, i_0$)}} &
				\multicolumn{2}{c}{\texttt{D($n_0, i_0$)}} \\
				\noalign{\smallskip}\hline\noalign{\smallskip}
				  $\eta_g$& Time (s)	& Gap (\%) &  Time (s)	& Gap (\%) &  Time (s)	& Gap (\%) & Time (s)	& Gap (\%) \\
				\noalign{\smallskip}\hline\noalign{\smallskip}  
		1	&	43.51	&	0.00	&	44.01	&	0.00	&	$\dagger$	&	45.17	&	$\dagger$&	69.59	\\
		1.1	&	27.88	&	0.00	&	26.88	&	0.00	&	$\dagger$	&	59.44	&	$\dagger$	&	69.33	\\
		1.2	&	26.63	&	0.00	&	26.84	&	0.00	&	$\dagger$	&	14.54	&	$\dagger$	&	69.51	\\
		1.3	&	22.19	&	0.00	&	30.55	&	0.00	&	$\dagger$	&	34.22	&	$\dagger$&	69.81	\\
		1.4	&	29.75	&	0.00	&	30.51	&	0.00	&	$\dagger$	&	6.91	&	$\dagger$	&	69.95	\\
		1.5	&	330.88	&	0.00	&	208.22	&	0.00	&	$\dagger$	&	2.58	&	$\dagger$&	71.69	\\
		1.6	&	$\dagger$	&	2.10	&	$\dagger$	&	2.09	&	$\dagger$	&	2.13	&	$\dagger$	&	71.43	\\
		1.7	&	$\dagger$	&	2.05	&	$\dagger$	&	3.84	&	$\dagger$	&	2.11	&	$\dagger$	&	71.73	\\
		1.8	&	$\dagger$	&	6.16	&	$\dagger$	&	7.80	&	$\dagger$	&	6.68	&	$\dagger$	&	71.86	\\
		1.9	&	$\dagger$	&	7.43	&	$\dagger$	&	7.62	&	$\dagger$	&	7.49	&	$\dagger$	&	71.80	\\
		2	&	$\dagger$	&	3.75	&	$\dagger$	&	3.81	&	$\dagger$	&	3.77	&	$\dagger$&	67.66	\\
		2.1	&	$\dagger$	&	5.04	&	$\dagger$	&	5.15	&	$\dagger$	&	5.05	&	$\dagger$	&	68.12	\\
		2.2	&	$\dagger$	&	5.01	&	$\dagger$	&	5.15	&	$\dagger$	&	5.01	&	$\dagger$	&	67.27	\\
		2.3	&	12.44	&	0.00	&	13.75	&	0.00	&	73.32	&	0.00	&	$\dagger$	&	67.84	\\
	\noalign{\smallskip}\hline
		\end{tabular}}
	\end{table}

	\begin{figure}[!t]
		\centering
		\subfloat[\texttt{D}, ($\eta_p,\eta_g$) = (1,1.2). ]{\includegraphics[width=0.35\textwidth]{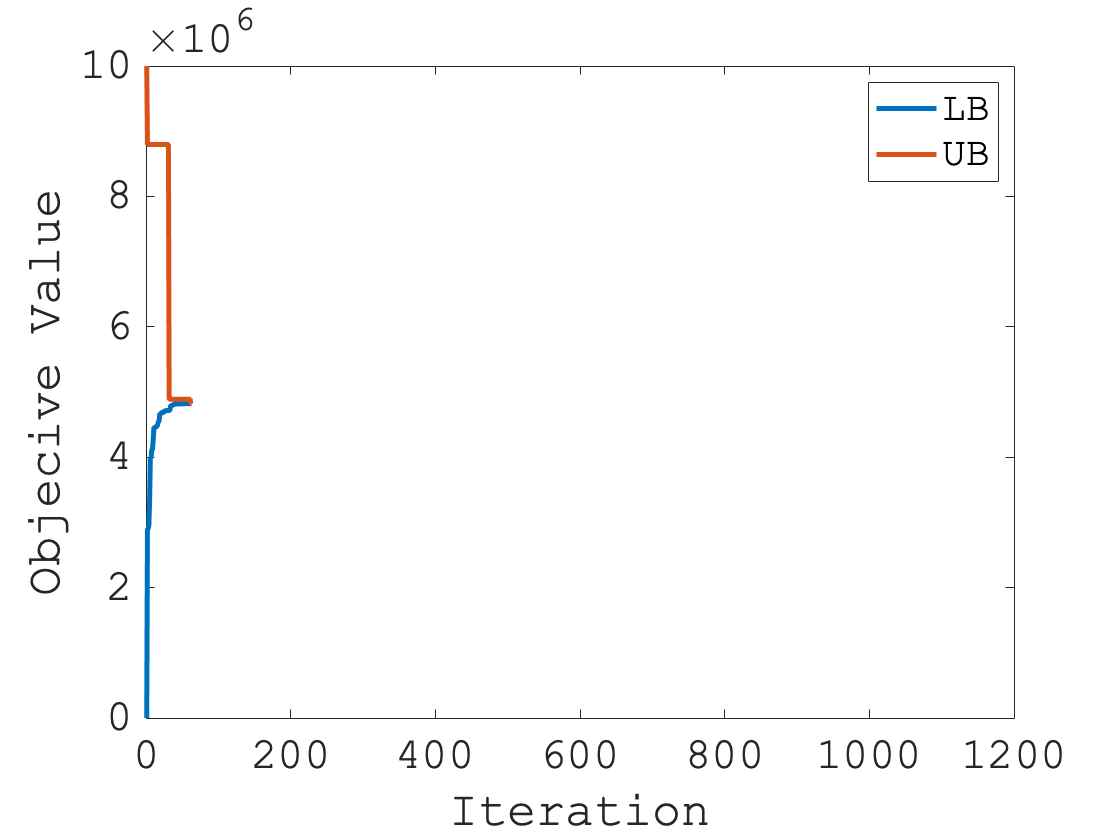}\label{fig:conv1D}}
		\subfloat[\texttt{B}, ($\eta_p,\eta_g$) = (1,1.2).]{\includegraphics[width=0.35\textwidth]{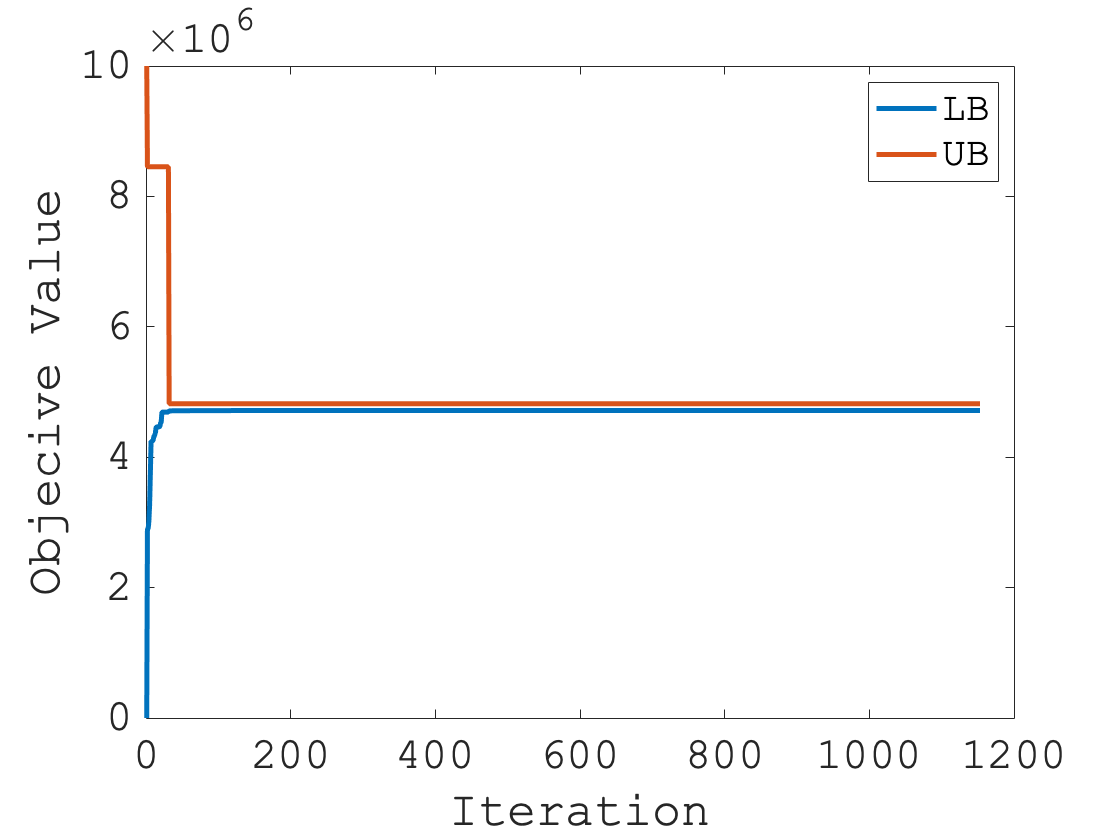}\label{fig:conv1B}}\\
		\subfloat[\texttt{D}, ($\eta_p,\eta_g$) = (1.6,1.8).]{\includegraphics[width=0.35\textwidth]{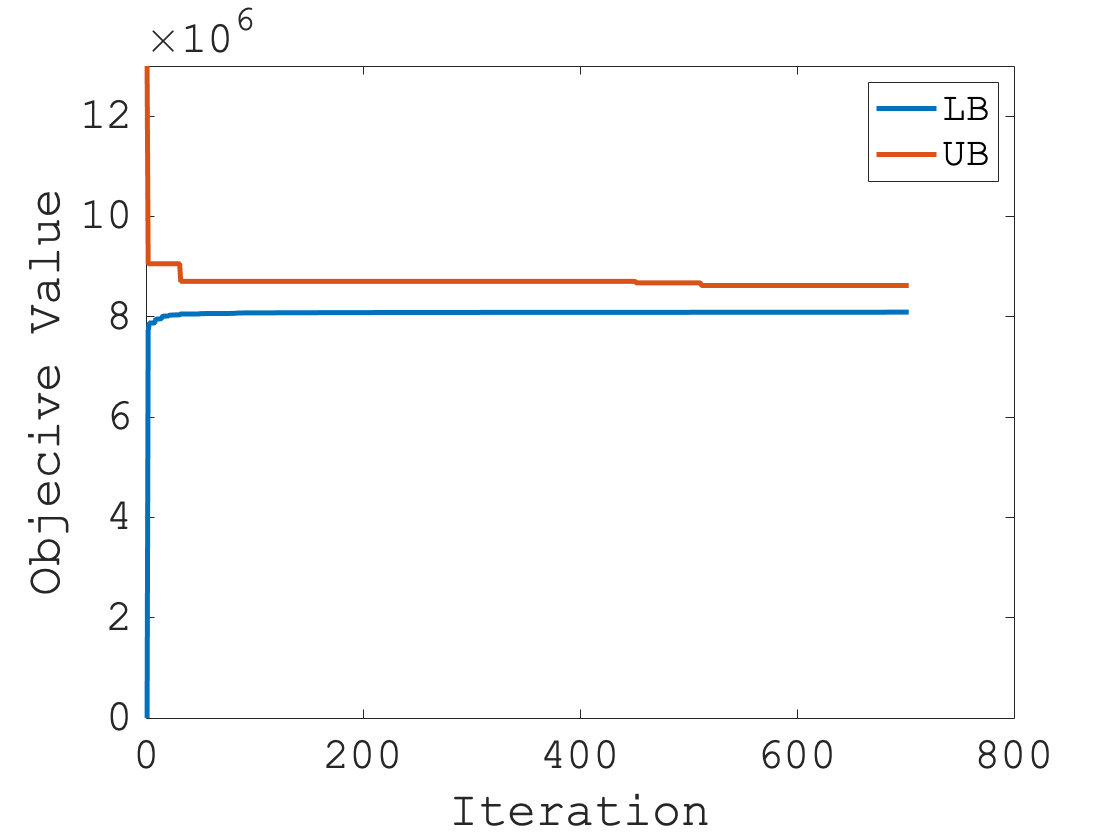}\label{fig:conv2D}}
		\subfloat[\texttt{B}, ($\eta_p,\eta_g$) = (1.6,1.8).]{\includegraphics[width=0.35\textwidth]{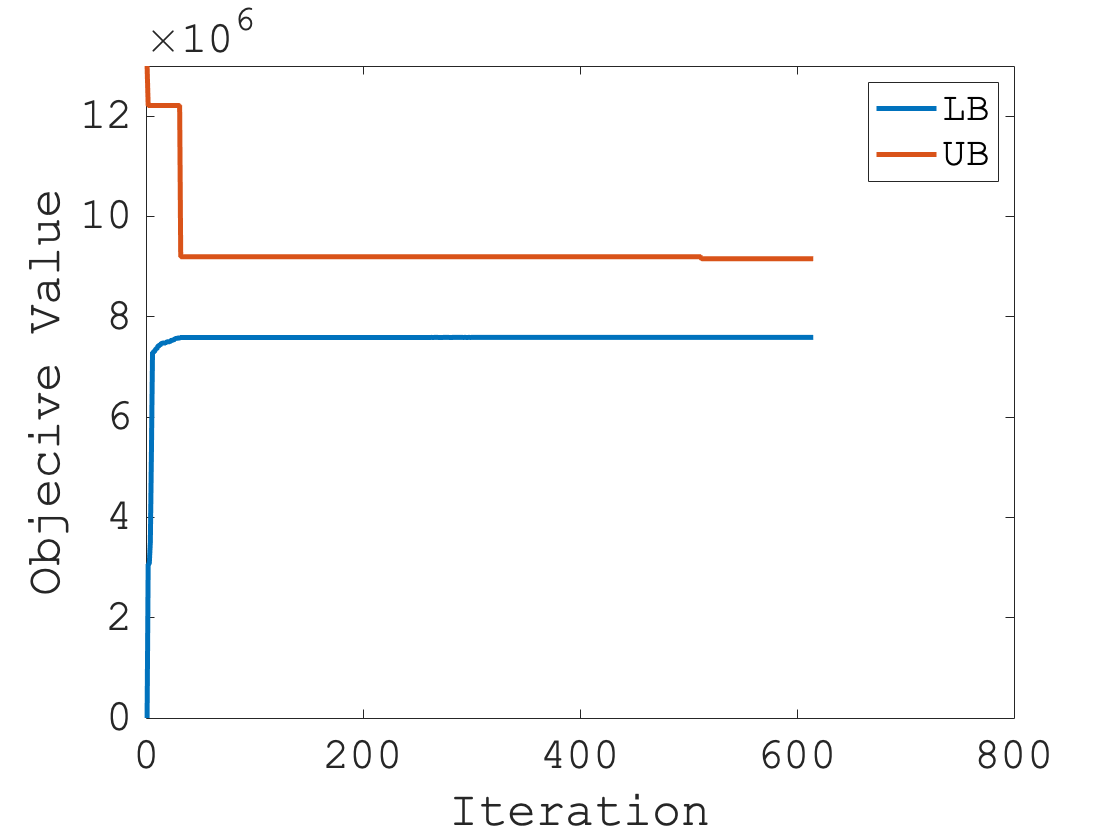}\label{fig:conv2B}}
		\caption{Convergence Behaviors of \texttt{D} and \texttt{B}.\label{fig:conv}}
		\end{figure}
\end{APPENDICES}

\end{document}